%
%

\magnification=1200

\font\headfont=cmr10 at 12 pt

\font\fpt=cmr10 at 14 pt



\def\AAA{\bf}
\def\BBB{\bf}

\overfullrule=0in

\def\boxit#1{\hbox{\vrule
 \vtop{%
  \vbox{\hrule\kern 2pt %
     \hbox{\kern 2pt #1\kern 2pt}}%
   \kern 2pt \hrule }%
  \vrule}}
  \def\mathqed{  \vrule width5pt height5pt depth0pt}

  \def\harr#1#2{\ \smash{\mathop{\hbox to .3in{\rightarrowfill}}\limits^{\scriptstyle#1}_{\scriptstyle#2}}\ }

 \def\GG{{{\bf G} \!\!\!\! {\rm l}}\ }

\def\ss{\subset}

\def\smfrac#1#2{\hbox{${#1\over #2}$}}

\def\log{{\rm log}}
\def\Hess{{\rm Hess}}

\def\tr{{\rm tr}}
\def\max{{\rm max}}

\def\det{{\rm det}}

\def\Sym{{\rm Sym}^2}

\def\rn{\bbr^n}

\def\Int{{\rm Int}}

\def\Symn{{\Sym(\rn)}}

\def\Theorem#1{\medskip\noindent {\bf THEOREM \bf #1.}}
\def\Prop#1{\medskip\noindent {\bf Proposition #1.}}
\def\Cor#1{\medskip\noindent {\bf Corollary #1.}}
\def\Lemma#1{\medskip\noindent {\bf Lemma #1.}}
\def\Remark#1{\medskip\noindent {\bf Remark #1.}}
\def\Note#1{\medskip\noindent {\bf Note #1.}}
\def\Def#1{\medskip\noindent {\bf Definition #1.}}

\def\Ex#1{\medskip\noindent {\bf Example \bf    #1.}}

\def\pf{\medskip\noindent {\bf Proof.}\ }
\def\qed{\hfill  $\vrule width5pt height5pt depth0pt$}

\def\cv{{\cal V}}   \def\cp{{\cal P}}

\def\cp{{\cal P}}
\def\cf{{\cal F}}

\def\vf{\varphi}

\def\wt{\widetilde}

\def\and{\qquad {\rm and} \qquad}

\def\ol{\overline}
\def\bbr{{\bf R}}\def\bbh{{\bf H}}
\def\bbc{{\bf C}}

\def\bbz{{\bf Z}}
\def\bbp{{\bf P}}

\def\a{\alpha}
\def\b{\beta}
\def\d{\delta}
\def\e{\epsilon}

\def\l{\lambda}

\def\s{\sigma}

\def\D{\Delta}

\def\O{\Omega}

\def\dbar{\ol{\partial}}

\def\bo{\partial \Omega}

\def\Symn{\Sym(\rn)}
 
\def\USC{{\rm USC}}

\def\cpt{\wt{\cp}}
\def\ft{\wt F}
\def\ob{\overline{\O}}

\def\HLDD{HL$_1$} 
\def\HLPUP{HL$_2$} 
\def\HLDDR{HL$_3$}
\def\HLREST{HL$_4$} 
\def\HLPCON{HL$_{5}$} 
\def\HLSURVEY{HL$_{6}$} 
\def\HLRS{HL$_{7}$}
\def\HLBP{HL$_{8}$} 
\def\HLAE{HL$_{9}$} 
\def\HLTangI{{HL$_{10}$}}
\def\HLTangII{{HL$_{11}$}}

\def\T{\Theta}

\def\ps{h}
\def\hub{\overline{h}}
\def\hlb{\underline{\psi}}
\def\Hz{H^0}

\def\AAA{1}
 \def\BBB{2}
 \def\BB{3}
 \def\CC{4}
 \def\DD{5}
  \def\AA{6}
  
\def\EE{7}
\def\FF{8}

\centerline
{
\fpt  THE DIRICHLET PROBLEM }
\medskip

\centerline{\fpt WITH PRESCRIBED INTERIOR SINGULARITIES }

\vskip.2in

\centerline{\fpt F. Reese Harvey and H. Blaine Lawson, Jr.$^*$}
\vglue .9cm
\smallbreak\footnote{}{ $ {} \sp{ *}{\rm Partially}$  supported by
the N.S.F.  and I.H.E.S.}

 
\centerline{\bf ABSTRACT} \medskip

\vskip .1in
  {{\parindent= .54in
\narrower  
In  this paper we solve  the nonlinear Dirichlet problem (uniquely) for functions with prescribed
asymptotic singularities at a finite number of points, and with
arbitrary continuous boundary data, on a domain in $\rn$.
The main results apply, in particular,  to subequations with a Riesz characteristic $p\geq2$.
It is shown  that, without requiring uniform ellipticity,  the Dirichlet problem can be 
solved uniquely
for arbitrary continuous boundary data with  singularities asymptotic to the 
Riesz kernel $\Theta_j K_p(x-x_j)$ where
$$
K_p(x) \ =\ \cases
{
- {1\over |x|^{p-2}} \qquad \  {\rm for}\ \ 2 < p< \infty,  \cr
\quad \log |x| \qquad\   {\rm if}\ \ p=2.
}
$$
at any prescribed finite set of points $\{x_1,...,x_k\}$  in the domain
and any finite set of positive real numbers $\Theta_1,..., \Theta_k$.  This  sharpens 
a previous result of the authors concerning the discreteness of high-density
sets of subsolutions. 

Uniqueness and   existence results are also established for finite-type singularities
such as $\Theta_j|x-x_j|^{2-p}$ for $1\leq p<2$.

The main results apply similarly with  prescribed
singularities asymptotic to the fundamental solutions of Armstrong-Sirakov-Smart
(in the uniformly elliptic case).

}}

\vskip.3in

\centerline{\bf TABLE OF CONTENTS} \bigskip

{{\parindent= .1in\narrower  \noindent

\qquad\qquad \AAA.  Introduction and Statement of  Some  Main Results.   \smallskip

\qquad\qquad \BBB. Riesz Kernels and the Riesz Characteristic.  \smallskip

\qquad\qquad \BB. Examples of Downward-Pointing $F$-Harmonics.  \smallskip

\qquad\qquad \CC. The Dirichlet Problem with Prescribed Densities.  \smallskip

\qquad\qquad \DD. Comparison.  \smallskip

\qquad\qquad \AA.  A Basic Construction and the Proof of Existence in the Polar Case.  \smallskip

\qquad\qquad \EE.  The Existence Proof in the Finite Case.\smallskip

\qquad\qquad \FF.  Prescribing Values at Singularities in the Finite Case.\smallskip

\qquad\qquad  Appendix A.   Asymptotic Equivalences and Tangent Flows.

}}

\vfill\eject


\centerline{\headfont \AAA.\ Introduction and Statement of Some Main Results.}
\bigskip

The aim of this paper is to study the Dirichlet problem for functions with prescribed
asymptotic singularities on a domain in $\rn$.  We shall adopt the notation and definitions
in our previous work [HL$_{1,3,6,10}$].


Throughout the paper $F\ss \Symn$ will denote a closed set which satisfies the weakest possible
ellipticity condition:

\medskip
\qquad
(F1). \ $A\in F$ and $P\geq0 \ \ \Rightarrow\ \ A+P \in F$,
\medskip
\noindent
so that solutions can be taken in the viscosity sense (cf. [C], [CIL], [CC]).
Said differently, $F$ is a constant coefficient, pure second-order subequation in $\rn$.

In addition we will always  assume that:

\medskip
\qquad
(F2). \ $F$ is a cone with vertex at the origin (i.e., $tF=F$ for $t>0$).
\medskip
\noindent
We will refer to an $F$ satisfying  (F1) and (F2) succinctly as a {\bf cone subequation}.

For certain existence results we shall also make the mild requirement that 
\medskip
\qquad
(F3). \  $-P_e \notin F$ \ for all unit vectors $e\in \rn$
\medskip\noindent
where $P_e$ denotes orthogonal projection onto the line $\bbr\cdot e$.
Said differently,  quadratic functions such as $u(x) = -x_1^2$ are not allowed to be 
$F$-subharmonic.
This property (F3) is equivalent to the following condition on the dual subequation  $\ft \equiv \sim(-F)$:
\medskip
\qquad
(B3). \  All smooth boundaries are strictly $\ft$-convex. 
\medskip\noindent
 See Proposition \AAA.7 below for the proof of this equivalence and further discussion.

Our asymptotic singularities will be prescribed by functions of the following type.

\def\pss{\psi}

\Def {\AAA.1} A {\bf downward-pointing singular $F$-subharmonic} at a point $x_0\in \rn$
is a continuous $[-\infty, \infty)$-valued $F$-subharmonic  function $\pss$ defined on a neighborhood
$U$ of $x_0$ with $\pss(x) > \pss(x_0)$ for $x\in  U-\{x_0\}$, such that either:
\medskip

{\bf (Polar Case):} \ \ \ \ \ $\pss(x_0) = -\infty$, \ \ \ or

\medskip

{\bf (Finite Case):} \ \ \ \ $\pss(x_0) > -\infty$, and $\pss$ has no test functions at $x_0$.

\medskip
\noindent
If in addition $\pss$ is $F$-harmonic  on $U-\{x_0\}$, then $\pss$ will be referred to as a 
 {\bf downward-pointing singular $F$-harmonic} at  $x_0$, and will usually be denoted by
 $h$ instead of $\pss$.

\medskip
\noindent
Note that   in the polar case $\pss$ also has no test functions at $x_0$.

\bigskip

\vfill\eject

The problem we want to address  is the following, which we will refer to as the {\bf (DPPS)}.
\medskip

\centerline
{\bf THE  DIRICHLET PROBLEM WITH PRESCRIBED SINGULARITIES.}
\smallskip
\noindent
Let $\O$ be a bounded domain in $\rn$, and fix  a finite number of points $x_1,...,x_k \in \O$.

\medskip
 
{\bf Boundary Data:} \hskip .15in
This consists of a function $\vf\in C(\bo)$ along with

\centerline{a downward-pointing
singular $F$-harmonic function $h_j$ at each point $x_j$.}

\medskip
\noindent
A function $H$ is a {\bf solution} of the corresponding (DPPS)   if 
\medskip

(1) \quad $H\in C(\ob - \{x_1,...,x_k\})$
\medskip

(2)  \quad  $H$ is $F$-harmonic on $\O - \{x_1,...,x_k\}$,
\medskip

(3)  \quad  $H\bigr|_{\bo} = \vf$,
\medskip

(4)  \quad  $H$ is {\bf asymptotically equivalent  to $h_j$ at $x_j$}.  By definition this means that:
\medskip

\qquad (4a) (In the Polar Case). There exists a constant $C>0$ such that for 

\qquad\qquad \ \ each  $j=1,...,k$
$$
h_j(x) -C \ \leq\ H(x)\ \leq\  h_j(x) +C  \qquad  {\rm near} \ x_j.
$$
\medskip

\qquad (4b) (In the Finite Case).  For each  $j=1,...,k$
$$
 \lim_{x\to x_j} {H(x) - H(x_j)  \over h_j(x) - h_j(x_j)}\ =\ 1.
$$

\Remark{\AAA.2}
One easily verifies that asymptotic equivalence is indeed an equivalence relation
(in fact, on the larger space of upper semi-continuous functions defined in a neighborhood of the point in question).
Each equivalence class  is invariant under the change of functions by an additive constant.   
We will denote {\sl  asymptotic equivalence}  
\medskip
\centerline{by $H \approx h$ in the polar case, \ \ and by $H\sim h$ in the finite case.}
\medskip

\medskip

The following result is an immediate consequence of
the comparison Theorems \DD.2 and \DD.4 proved in Section \DD.
\medskip

\Theorem{\AAA.3. (Uniqueness)} 
{\sl  For any  cone subequation $F$  
there is at most one function
$H$ with properties (1) through (4).}

\bigskip

An existence construction is presented in Section  \AA\ and then completed in the polar case
by proving the following theorem.  
Our existence results in the finite case are stated and proved in Section \EE.
\medskip

\Theorem {\AAA.4. (Existence in the Polar Case)}  {\sl
Suppose that $F$ is a cone subequation satisfying Condition (F3) = (B3),
and that  the boundary of $\O$  is smooth and strictly  $F$-convex.
Assume that  each prescribed singularity $h_j$ at $x_j$ is a downward-pointing singular $F$-harmonic 
of polar type. Finally assume
\medskip

\noindent
{\bf Hypothesis (H):} 
There exists a continuous $F$-subharmonic function $\pss$ on $\ob$,
finite except at $x_1, ... , x_k\in \O$,
with $\pss\approx h_j$ at $x_j$ for each $j$.

\medskip
\noindent
 Then the   Dirichlet  Problem with Prescribed Singularities described above has a  solution.
 Moreover, it is uniquely determined as the Perron function
 $$
 H(x) \ =\ \sup_{v\in \cf} v(x) 
 $$
for the family 
\medskip
 
\hskip.8in
 $\cf \equiv \{v\in \USC(\ob) : v $ is $F$-subharmonic on $\O$, $v\leq \vf$ on $\bo$, 
 
\hskip 1.5in and $v- h_j$ is bounded above near each    $x_j\}$.
}

\bigskip

By Remark 11.13 in  [\HLDDR]   this Theorem extends to domains which are  finite intersections of domains with smooth strictly $F$-convex boundaries. 

Typically the functions $h_j$ arise from a single global $F$-subharmonic function
$h$ on $\rn$ which is $F$-harmonic outside the origin and has a downward-pointing
singularity at 0. Such an $h$ will be called a {\sl generalized fundamental solution}
for the subequation $F$.  It determines the  data in the (DPPS) by taking
$h_j(x) = \T_j h(x-x_j)$ with constants $\T_j>0$ for $j=1, ... , k$.  For Hypothesis (H) consider the function
$$
\psi(x) \ =\ \sum_{j=1}^k \T_j h(x-x_j).
\eqno{(\AAA.1)}
$$
In the polar case (where $h(0)=-\infty$) the condition $\psi \approx h_j$ at $x_j$ is automatic
since $h$ is continuous outside the origin.  Thus the hypothesis (H) reduces to 
\medskip

\noindent
{\bf Hypothesis (H)$'$:}  \ \ $ \psi(x) \ = \ \sum_j \T_j h(x-x_j)$ is $F$-subharmonic on $\rn$.

\medskip

There are two cases where this hypothesis is easily satisfied.
\medskip
\noindent
{\sl A Single Point Singularity at $x_0\in \O$:} Then $\psi(x) \equiv \T h(x-x_0)$ ($\T>0)$ is $F$-subharmonic
on $\rn$ -- in fact, $F$-harmonic on $\rn-\{x_0\}$.

\medskip
\noindent
{\sl The subequation $F$ is convex:} Then $\psi$ is $F$-subharmonic on $\rn$ because sums and positive multiples
of $F$-subharmonic functions are also $F$-subharmonic.

\medskip
This yields two special cases where the (DPPS) is uniquely solvable.

\Cor{\AAA.5}
{\sl Let $F$ and $\O$ be as in Theorem \AAA.4, and 
suppose $h$ is a generalized fundamental solution for $F$.
\smallskip

(a)\ \ For each $x_0\in \O$ and $\T>0$ there exists a unique solution to the (DPPS) 
having boundary values $\vf$ on $\bo$ and asymptotic to $\T h(x-x_0)$ at $x_0$.

\smallskip
(b) \ \ If, in addition, $F$ is convex, then the multi-pole (DPPS)  with $\psi$ as in (\AAA.1), has a unique solution.
}
\medskip

In the case of just one point singularity with $\T = 1$ and the outer boundary function $\vf\equiv 0$, the 
solution provided by Corollary \AAA.5(a) will be denoted $G_\O(x; \, x_0, h)$ and referred to as the
{\bf nonlinear Green's function} for the subequation $F$ on the domain $\O$
with asymptotic  singularity determined by  the generalized fundamental solution $h$.

In the multi-pole case where $F$ is required to be convex, we again take $\vf\equiv 0$.  Then the
solution to the (DPPS) given by Corollary \AAA.5(b) will be called the {\bf multi-pole nonlinear Green's function}
and denoted by $G_\O(x;\, x_1, ... ,x_k; \, \T_1, ... , \T_k; \, h)$.  

These functions extend the classical  pluri-complex Green's function
with logarithmic singularities (cf.  [Lem], [K$_*$], [Le], [Z])   where $h(x) = \log |x|$ on $\bbr^{2n} = \bbc^n$.
(See Theorem \BB.5 ff. for further discussion.)

In Section \EE \ we establish existence results in the finite case.  The strongest result,
Theorem \EE.4, applies when there is just one singularity. 
Here the asymptotic type of the singularity can be prescribed along with the outer boundary values.  
One begins with a function $h\in C(\ob)$ which is $F$-harmonic
on $\O-\{x_0\}$ and has a downward-pointing singularity  of finite type at $x_0$.
One then says that a function $H\in C(\O)$ has {\bf $h$-density} $\Theta\geq0$ if
$$
\lim_{x\to x_0} {H(x)-H(x_0)  \over h(x)-h(x_0)}\ =\ \Theta.
$$
This is equivalent to saying that $H \sim \Theta h$ at $x_0$.
Theorem \EE.4 asserts that under certain assumptions on $F$, the Dirichlet problem
can be solved for a harmonic with any prescribed $h$-density $\Theta\geq0$ at $x_0$ and $\vf\in C(\O)$.
It is then shown that the hypotheses in Theorem \EE.4 are satisfied for two large and
important classes of subequations. The first is the class of O$(n)$-invariant subequations
whose Riesz-characteristic $p$ satisfies $1\leq p <2$.  Here $h(x) = |x-x_0|^{2-p}$.
The second class consists of the uniformly elliptic  subequations where $h$
is taken to be the downward-pointing fundamental solution of Armstrong, Sirakov and Smart. 

There is another existence question which is meaningful only in the finite case.
As before we fix a boundary function $\vf\in C(\bo)$ and points $\{x_1,...,x_k\}\ss \O$. 
However, rather  than prescribing  densities, 
we instead prescribe the values $v_j$ for $H$ at each point $x_j$.  This is 
the Dirichlet Problem on the Punctured Domain $\O-\{x_1,...,x_k\}$,
 where we look for a function $H\in C(\ob)$ which is
$F$-harmonic on   $\O-\{x_1,...,x_k\}$ with $H\bigr|_{\bo}=\vf$ and $H(x_j)=v_j$
for each $j$.    By comparison on  $\O-\{x_1,...,x_k\}$,
if a solution $H$ exists, it is unique (see for example, [\HLAE, Thm. 6.2]).
This leaves the important problem of exactly determining the set  $\cv\!{\it al}$ of values $v=(v_1,...,v_k)$ 
for which a solution $H$ exists.
In general  $\cv\!{\it al}\ss\bbr^k$ is a proper subset which depends on the given data.
In Section \FF \ we establish the existence of a large, and explicitly described, subset
$\cv\ss \cv\!{\it al}$.  For the case where $\vf =0$,
it is shown that this set $\cv$ 
 is a convex cone with vertex at the origin and non-empty interior contained in the  negative ``octant''  $\bbr_-^k$.
When $\vf=0$ and there is only one point ($k=1$), the ``value problem'', namely, that of determining
$\cv\!{\it al}$, is solved completely:   $\cv\!{\it al} = \{v\leq0\} = \cv$ for any choice of point $x_1\in \O$.
(See Proposition \FF.5.)

The main theorem in Section \FF\ (Theorem \FF.1) actually enables one to prescribe not only certain
values of $H$ at the points $x_j$ but also to prescribe the tangent to $H$ at $x_j$
up to a positive multiple $\geq 1$.  The key hypothesis in Theorem \FF.1 is that there must exist on
$\ob$ an $F$-subharmonic function $h$ which is $\leq \vf$ on $\bo$ and asymptotic
to the given tangent at each $x_j$.  For convex subequations this is often easily done
and one obtains the large subset $\cv\ss \cv\!{\it al}$ discussed above.

\vfill\eject

\Remark {\AAA.6. (The Boundary Convexity Hypothesis)}
The hypothesis that $\bo$ is strictly $F$-convex in Theorem \AAA.4 is necessary
for existence  in the finite cases as well as
the polar case.  Of course there are many domains with smooth boundary
where this is true.  For example, if $\bo$ is strictly convex, then $\bo$ is also
strictly  $F$-convex for any cone subequation $F$  because $\cp\ss F$.
On the other hand, there are many cone subequations  $F$  with the property that every smooth  boundary is
 strictly $F$-convex, allowing $\bo$ to be arbitrary.  This is true  if and only if $P_e\in \Int F$ for all $|e|=1$ by the following Proposition,  applied to $\ft$ instead of $F$.

\medskip
\noindent
{\bf Proposition \AAA.7. (Concerning Condition (F3)).}
{\sl
For a cone subequation $F$ the conditions:
\medskip

(F3) \ \ $-P_e \notin F$  for any unit vector $e$, and

\medskip

(B3)\ \ All boundaries are strictly $\ft$-convex

\medskip
 \noindent
 are equivalent.    Furthermore,  if $F$ is  a convex subequation, then (F3) and (B3) are  also equivalent to:
 \medskip

(F3)$'$\ \  The subequation $F$ is {\sl complete}, i.e., it cannot be defined using the variables in a proper 
linear subspace of $\rn$.
\medskip

Finally, if $F$ is an invariant cone subequation as  in Lemma \BB.4, then these conditions are equivalent to
\medskip

(F3)$''$  \ \ $F$ has  finite Riesz characteristic $p$ (see Definition \BB.2).
\medskip
}


We shall use condition (B3), while condition (F3) provides the simplest 
test for this assumption.  Condition (F3)$'$ is automatic unless the subequation is so
degenerate that it  does not involve all the variables in $\rn$.

\pf
For cone subequations the equivalence of    (F3)  and  (B3)  was established in Section 5 of the earlier paper
[\HLDD], but was not stated explicitly.  The argument goes as follows.  First, by Lemma 5.3 (ii)$'$ in [\HLDD] condition 
(B3) is equivalent to 
\medskip

(F3)$'''$ \ \ For all $B\in \Symn$ and all unit vectors $e$\medskip
\centerline
{
$B+tP_e \in \Int \ft$ for all $t\geq $ some $t_0$.
}
\medskip\noindent
(This can be considered to be the definition of strict $\ft$-convexity.)
Second, by the  Elementary Property (5) in [\HLDD, \S 3] with $B'\equiv 0$, this is equivalent to
\medskip

(F3)$^{(iv)}$ \ \ $P_e \in \Int \ft$ for all unit vectors $e$.
\medskip
\noindent
Finally, by definition we have $-\Int \ft =\ \sim F$, so that (F3)$^{(iv)}$ $\iff$ (F3).

The equivalence of the two structural conditions (F3) and (F3)$'$, for $F$  convex, 
was  established in Proposition 3.6 of [\HLPUP].

For the equivalence  of (F3) and (F3)$''$ see Lemma \BB.4.   \qed

\vfill\eject


\centerline{\headfont
\BBB. Examples of Downward-Pointing F-Harmonics.
}

\Ex{\BBB.1. (Riesz Kernels)} Perhaps the most important such examples are the
classical Reisz kernels $K_p$ (with $p\geq1$) which are defined and discussed
in the next section.  Each has a downward-pointing singularity.
They play a fundamental role in standard potential theory (cf.\ [L]).
  Moreover, each $K_p$ is  actually a punctured harmonic for a large family of subequations -- those of {\sl Riesz characteristic} $p$. For such subequations $F$ the Riesz kernels are  central to the study of tangents and densities in the associated $F$-potential theories [\HLTangI], [\HLTangII]. 
  In particular, they often arise  as the unique tangent to any $F$-subharmonic function.

     A large and important class of  subequations with characteristic-$p$, which are convex cones
     but not uniformly elliptic, come from:
    
    \smallskip
    \noindent  
    {\bf Geometrically Defined Subequations:}
These are the convex cone subequations $\cp(\GG)$ determined by a closed subset
$\GG\ss G(p,\rn)$, of the Grassmannian of $p$-dimensional subspaces of $\rn$,
by the requirement that
$$
A\in F \qquad\iff\qquad \tr\left( A\bigr|_W  \right) \ \geq\ 0 \qquad \forall\, W\in\GG.
$$
The following {\bf Fullness Condition on $\GG$}:
$$
{\rm Each \ unit\  vector \ \ } e\in \rn\ \ {\rm is\ contained\ in \ a\  subspace\ \ } W\in \GG,
\eqno{(\BBB.1)}
$$
is equivalent to the $p^{\rm th}$ Riesz kernel being a punctured $\cp(\GG)$-harmonic on $\rn-\{0\}$.
Moreover,  if $\GG$ satisfies (\BBB.1),
then  Condition (F3) = (B3) also holds for $\cp(\GG)$. 
 
 These examples contain the subequations naturally associated
to many calibrations.  They also include  the Lagrangian subequations on $\bbc^n$.
If $\GG=G(p,\rn)$, the resulting subequation $\cp(\GG)$ (denoted $\cp_p$ in Example \BB.6(1)) is basic in geometry
(cf. [Wu], [Sha] and [\HLPCON,] for example).
See Section 4 of  [\HLTangI]   and Appendix A in [\HLTangII] for many more examples.

Perhaps it deserves mentioning here that the potential theory associated with the subequation
$\cp(\GG)$ is more appropriately called the {\sl $\GG$-pluripotential theory} because
of the fact that: $u$ is $\cp(\GG)$-subharmonic $\ \Leftrightarrow\ $ $u\bigr|_W$ is $\D$-subharmonic
for every affine $W\in\GG$ (proved in [\HLREST]).

\Ex{\BBB.2. (Homogeneous Singularities)}
 Suppose that the function $\pss$ in Definition \AAA.1 is homogeneous.
This means (assuming $x_0=0$ for simplicity) that  
$$
\pss(x) = |x|^\a \pss \left({x\over |x|}\right) \qquad {\rm for\ some\ \ } \a\neq 0
$$
 or, in the case ``$\a=0$'', $\pss(x) = \pss({x\over |x|}) + \T \log|x|$ 
where $\sup_{|x|=1} \pss = 0$ and $\T > 0$.  
Note that the Riesz kernel $K_p$ has homogeneity $\a = 2-p$.

Suppose that $\pss$ is a homogeneous $F$-subharmonic on $\rn-\{0\}$.
Then $\pss$ is downward-pointing (Definition \AAA.1) as follows.

\medskip
\noindent
{\bf Case: $0<\a \leq 1$.}
Then $\pss$ has a strict minimum at $x_0=0$ if and only if $\pss(x)>0$ for $x\neq0$, in which case $\pss(0)=0$
and $c|x|^\a \leq \pss(x)$  with $c=\inf_{\partial B} \psi$.   The condition of having no test functions
at 0  follows easily from this inequality since  $0<\a\leq 1$.  Thus $\pss$ has a  downward-pointing
singularity at 0 $\iff$ $\pss(x)>0$ on $|x|=1$.

\medskip
\noindent
{\bf Case: $\a<0$.}  Then $\pss$ has a strict minimum at $x_0=0$ if and only if $\pss(x)<0$ for $x\neq0$, in which
case $\pss(0)=-\infty$.  Thus $\pss$ has a  downward-pointing
singularity at 0 $\iff$ $\pss(x)<0$ on $|x|=1$.

\medskip
\noindent
{\bf Case: $\a=0$.} For any $\pss({x\over |x|})$ continuous we have $\pss(0)=-\infty$, while $\pss(x)$
is finite for $x\neq0$.  Thus $\pss$ always has a  downward-pointing
singularity at 0.

\Ex{\BBB.3. (The Special Cases $\cp$, $\cp^\bbc$ and $\cp^\bbh$)}  
These are the geometrically defined subequations obtained by taking $\GG$ to be $G(1,\rn)$,
$G^\bbc(1,\bbc^n) \ss G(2, \bbr^{2n})$ and 
$G^\bbh(1,\bbh^n) \ss G(4, \bbr^{4n})$ respectively.
Here there are huge families of downward-pointing $F$-harmonic functions $h$. Those which
are homogeneous can be classified as follows.  (The differential inequalities are in the viscosity sense.)

\medskip
\noindent
{\bf The $\cp$ Case:}
$$
\eqalign
{
h(x) \ =\ &|x| g\left({x\over |x|}\right) \qquad{\rm with}\ \ g\in C(S^{n-1}) \qquad {\rm satisfying}  \cr
&\Hess_{S^{n-1}} g + g I \ \geq\ 0, \quad{\rm and}\ \ \inf_{S^{n-1}} g \  > \ 0
}
$$

\medskip
\noindent
{\bf The $\cp^\bbc$ Case:}
$$
\eqalign
{
h(x) \ =\ & g\left( [x]\right) + \T \log |x|\qquad{\rm with}\ \ g\in C(\bbp(\bbc^n))  \ \ {\rm and}\ \ \T>0 \quad {\rm satisfying}  \cr
& i\partial\dbar g + \T \omega  \ \geq\ 0, \quad{\rm and}\ \ \sup_{\bbp(\bbc^n)} g \ =\ 0  \cr
\ {\rm (equivalently} \ \ \  &\Hess_{\bbp(\bbc^n)}^\bbc g + \T I \ \geq\ 0, \quad{\rm and}\ \ \sup_{\bbp(\bbc^n)} g \ =\ 0).
}
$$

\medskip
\noindent
{\bf The $\cp^\bbh$ Case:}
$$
\eqalign
{
h(x) \ =\ & {1\over |x|^2} g\left( [x]\right) \qquad{\rm with}\ \ g\in C(\bbp(\bbh^n)) \qquad {\rm satisfying}  \cr
&\Hess_{\bbp(\bbh^n)}^\bbh g -  2g I \ \geq\ 0, \quad{\rm and}\ \ \sup_{\bbp(\bbc^n)} g \ < \ 0
}
$$
See Section 5 in [\HLTangII] for details and proofs.

\Ex{\BBB.4. (Bedford-Taylor Singularities)}  
Downward-pointing $F$-harmonics need not be homogeneous.
The  $\cp^\bbc$-harmonic (i.e.,  maximal plurisubharmonic) functions  on $\bbc^n-\{0\}$ given by
$$
h(z,w)\ =\ \log(|z|^2 + |w|^4)
$$
and each of  its unitary rotates, define a family of distinct asymptotic equivalence 
classes at the origin, which are  different from the basic punctured harmonic $\log(|z|^2+|w|^2)$.
There are also the $\cp^\bbc$-harmonic functions $\log(|z|^\a+|w|^\b)$ for $\a,\b>0$.
These examples go back to Bedford and Taylor  [BT].

\Ex{\BBB.5. (Armstrong-Sirakov-Smart Singularities)} In the very interesting paper [AS$_1$]
the authors consider  pure second-order cone subequations  exactly as in this paper
 but with the additional hypothesis of 
uniform ellipticity  (where  (F3)  is automatic). Under this hypothesis they
 establish the existence of a canonical
{\sl fundamental solution} $\Phi$ which  is $F$-harmonic on $\rn-\{0\}$,  
 has the homogeneity property $\Phi(tx) = t^{2-p}\Phi(x), \forall\, t>0$ for some $1 < p<\infty$,
 and has a downward-pointing singularity at 0.
 In stark contrast to the three special cases $\cp, \cp^\bbc$ and $\cp^\bbh$ in Example \BBB.3
above, in this uniformly elliptic case there is precisely one asymptotic equivalence class
of downward-pointing $F$-harmonic, namely $\Phi$,  up to a positive scale  (see  [AS$_1$]).

Now start with any cone subequation $F\ss\Symn$ and consider the family 
$F(\d), \d>0$ of elliptic regularizations consisting of
uniformly elliptic cone subequations converging to $F$ as $\d\to0$ (See Appendix A in [\HLTangII]).
The paper [AS$_1$] applies to each $F(\d)$ presenting us with a large set of functions $h_\d \equiv \Phi$,
 other than the Riesz kernels, which 
satisfy the conditions of Definition \AAA.1.

Now suppose that $F$ is convex, so that for small $\d$, $F(\d)$ is also convex.  
The smallest possible subequation $F_\d$ for which the function  $h_\d = \Phi$
is an entire {\sl sub}harmonic is constructed by taking the closed convex cone on 
$\{D^2_x \Phi: |x|=1\}   + \cp$.  (See Example \BB.6 where this is carried out in the cases
$F=\cp, \cp^\bbc$ and $\cp^\bbh$.)
If $\Phi$ is also a harmonic for this smallest $F_\d$ (outside the origin),
then $\Phi$ is a downward-pointing harmonic for all cone subequations $F$
such that $F_\d \ss  F \ss F(\d)$.  For each such $F$, Theorems \AAA.3 and \AAA.4 apply,
as well as Theorem \EE.1 in the finite case.

Other examples of this phenomenon can be constructed by applying a non-orthogonal linear
transformation to any of the many convex subequations studied in [\HLTangI] and  [\HLTangII].


\vskip.3in


\centerline{\headfont
\BB.  Riesz Kernels and the Riesz Characteristic.
}
\medskip
 
Of particular importance for the study of isolated singularities of subsolutions are the
classical {\sl Riesz kernels}:
$$
K_p(x) \ =\ \cases
{
-{1\over p-2} {1\over |x|^{p-2}} \qquad \  {\rm for}\ \ 1\leq p< \infty, \ p\neq 2 \cr \cr
\quad \log |x| \qquad \qquad  {\rm if}\ \ p=2.
}
$$
One sees easily  that 
$$
D_x K_p \ =\ {x\over |x|^p} 
\and
D^2_x K_p \ =\ {1\over |x|^p} \bigl( P_{x^\perp} -(p-1) P_x\bigr)
$$
where $ P_{x^\perp} $ and $P_x$ denote orthogonal projection onto the hyperplane
perpendicular to $x$  and the line through $x$ respectively.  From this one sees
the following.

\Prop{\BB.1} {\sl
Given a cone subequation $F$, the $p^{\rm th}$ Riesz kernel}
$$
K_p \ \ {\sl is\ \ } F\,{\sl  harmonic \ on\ \ } \rn-\{0\} 
\qquad\iff\qquad 
P_{e^\perp} -(p-1)P_e \in \partial F \ \ \ \forall\, |e|=1.
\eqno{(\BB.1)}
$$

\Def{\BB.2} A (not necessarily convex) cone subequation $F\ss\Symn$ is said to have a
{\bf finite Riesz characteristic $p_F = p$} if 
$$
 P_{e^\perp} -(p-1) P_e \ \in\ \partial F \qquad{\rm for\ all\ unit\ vectors\ \ } e\in \rn.
$$

\medskip

 If $F$ has finite Riesz characteristic $p$, then 
$1\leq p$ since $\Int\cp\ss\Int F$.  In addition,
the Riesz kernel $K_p$
and  each of its translates is {\bf $F$-harmonic} outside its singularity.
Thus 
$$
h(x) \equiv K_p(x) \ \ {\rm  is \ a  \ downward\ pointing,\ singular \ } F\, {\rm harmonic\ on\ } \rn-\{0\}
\eqno{(\BB.2)}
$$
for any cone subequation with Riesz characteristic $p$.
 This is perhaps the most important
example satisfying the conditions in Definition \AAA.1. It is discussed
in detail after a few additional comments on the polar case $p\geq2$.

\Remark{\BB.3}  It is natural to extend the definition of Riesz characteristic.  We say that $F$ 
{\bf has Riesz characteristic} $p_F=\infty$ if
$$
-P_e \in \partial F, \ \ {\rm or\ equivalently}\ \ -P_e \in F, \quad \forall\, |e|=1.
\eqno{(\BB.3)}
$$
The equivalence is because $-P_e \notin \Int F$ for a cone subequation $F$
unless $F=\Symn$.

\medskip


Of course, there are cone subequations which do not have a Riesz characteristic.  However,
as seen in the next lemma,
for ``invariant''  cone subequations the Riesz characteristic is very easy to compute, 
and condition (F3) holds $\iff$ $p$ is finite.

\Lemma{\BB.4} {\sl
Suppose that $F$ is invariant under a subgroup $H\ss {\rm O}(n)$ which acts
transitively on the sphere $S^{n-1}\ss\rn$. Then
$$
p_F \ =\ \sup\{q :  P_{e^\perp} -(q-1) P_e \ \in\ F \}
\eqno{(\BB.4)}
$$
for some (and therefore any) unit vector $e\in \rn$.
In particular, $p_F =\infty \iff $ (F3) fails for all $|e|=1$.
}
\medskip

The proof is straightforward.  Note that in appropriate coordinates $P_{e^\perp} -(q-1) P_e$ is  diagonal matrix
with one eigenvalue $-(q-1)$ and the remaining entries  all equal to 1.

The principal Corollary \AAA.5 has as its most important special case the following.

\Theorem{\BB.5. (The Polar Case)} 
{\sl
Suppose $F$ is a cone subequation with  finite Riesz characteristic $p\geq 2$ and property (F3), and
let $\O$ be a domain with a smooth boundary which is strictly $F$-convex.
\smallskip

(a)\ \ For each $\vf\in C(\bo)$, $x_0\in \O$, and $\T>0$, there exists a unique solution $H$ to the (DPPS)
having boundary values $\vf$ on $\bo$ and asymptotic singularity $H(x)\approx \T K_p(x-x_0)$ at $x_0$.
\smallskip

(b)\ \ If, in addition, $F$ is convex, then the multi-pole (DPPS) with asymptotics  $\T_jK_p(x-x_j)$ has a 
unique solution $H$.
\smallskip

Moreover, in either of these cases, if $\vf\equiv0$, this provides the existence and uniqueness of a
nonlinear Green's function $G_\O$ with
$$
G_\O(x) \ \approx \ \T K_p(x-x_0)\qquad {\rm at\ \ } x_0
$$
in the case of a single pole, and, provided that $F$ is convex, 
$$
G_\O(x) \ \approx \ \T_j K_p(x-x_j)\qquad {\rm at\ each \ \ } x_j
$$
in the multi-pole case.
}
\medskip

The existence and uniqueness of the multi-pole Green's function  for the subequation $F=\cp^\bbc$ 
 was proved by Lelong in 1989 ([Le]).
This built on previous work for single-point Green's functions (Lempert [Lem] and Klimek [K$_*$]).
An even more general version was established by A. Zeriahi  [Z].

Theorem \BB.5  includes all invariant cone subequations whose Riesz characteristic  is finite
(by Lemma \BB.4).

There are many more subequations $F$ which have a finite Riesz characteristic than one
might at first imagine.  We start by mentioning four extreme examples of characteristic $p$ subequations.
The first explains our choice of normalization in the definitions of $K_p$ and characteristic $p$.  When $p$ is an integer, this example coincides
with the geometric subequation $\cp(G(p, \rn))$ discussed in Section 2.

\medskip
\noindent
{\bf Examples. \BB.6.}   
\medskip
\noindent
(1) \ \ $\cp_p  = \{A : \l_1(A)+\cdots+ \l_{[p]}(A) +(p-[p])\l_{[p]+1}(A) \geq 0\}$  where 
$ \l_1(A)\leq  \l_2(A)\leq\cdots $ are the ordered eigenvalues of $A$ and $1\leq p\leq n$.

\medskip
\noindent
(2) \ \ $\cp(\d_p) = 
\left    \{A : A+ \smfrac  {\d_p}{n}  \tr(A)  I \ \geq\ 0    \right  \}$, where $\d_p =  {n(p-1) \over n-p}$.

\medskip
\noindent
(3) $\cp_p^{\rm min/max} = \{A : \l_{\rm min}(A) + (p-1) \l_{\rm max}(A) \geq 0\}$.

\medskip
\noindent
(4) $\cp_p^{\rm min/2} = \{A : \l_{\rm min}(A) + (p-1) \l_{2}(A) \geq 0\}$.

\medskip

Note that $\cp_p^{\rm min/2} \ss\cp_p \ss \cp(\d_p)\ss\cp_p^{\rm min/max}$
and that each subequation $F$ with $\cp_p^{\rm min/2} \ss F \ss\cp_p^{\rm min/max}$
has finite Riesz characteristic $p$.  
It is somewhat surprising that, under a mild restriction, there exist both a ``largest'' and a 
``smallest'' characteristic $p$ subequation.
More precisely,  {\sl every invariant (as in Lemma \BB.4)
cone subequation $F$
with finite Riesz characteristic $p$ satisfies 
$$
\cp_p^{\rm min/2}  \ \ss\   F  \ \ss\  \cp_p^{\rm min/max}
\eqno{(\BB.5)}
$$
and if $F$ is O$(n)$-invariant and  convex, then}
$$
\cp_p \ \ss\   F  \ \ss\  \cp(\d_p).
\eqno{(\BB.6)}
$$
This is proved in Appendix A of [\HLTangII] where many more examples
of characteristic $p$ subequations are given. These include subequations
of Monge-Amp\`ere type arising from G\aa rding operators. Among these is the following
{\sl Hessian equation}, which has been studied by Trudinger-Wang [TW$_*$], Labutin [La$_*$] and others.
We have drawn heavily from [La$_3$] in this paper.
\medskip
\noindent
(5)  $\Sigma_k=\{A : \s_1(A)\geq0,...,\s_k(A)\geq 0\}$, where  $p_F ={n  \over k}$.
\medskip
\noindent
One also has
\medskip
\noindent
(6) $F= \{A : \tr(A^q)\geq0\}$ where $p_F =1 + (n-1)^{1\over q}$ for $q\in \bbz$ odd.
(This $F$ is not convex.)

\medskip
\noindent
Suppose $F$ is an O$(n)$-invariant subequation with $p_F =p$, and let 
$F(\bbc), F(\bbh)$ be the complex and quaternionic analogues given by the
same conditions on the eigenvalues of their hermitian symmetric components.
Then $p_{F(\bbc)}=2p_F$ and $p_{F(\bbh)}=4p_F$.
These examples  contain the complex and quaternionic Monge-Amp\`ere equations 
  [BT], [A$_*$], [AV], as well as the  
   complex and quaternionic  hessian equations. There are also the complex
   (as well as quaternionic) analogues:
   $$
   \cp_p^{\rm min/2}(\bbc)  \ \ss\  \cp_p(\bbc)   \ \ss\   \cp(\d_p)(\bbc) \ \ss\  \cp_p^{\rm min/max}(\bbc).
   \eqno{(\BB.7)}
$$
with
$$
\cp_p^{\rm min/2}(\bbc)  \ \ss\   F(\bbc)  \ \ss\  \cp_p^{\rm min/max}(\bbc)
\and
\cp_p^{\rm min/2}(\bbc)  \ \ss\   F(\bbc)  \ \ss\  \cp_p^{\rm min/max}(\bbc)
$$
as in (\BB.5) and (\BB.6).


\vskip .3in

\centerline{ \headfont \CC.  The Dirichlet Problem with Prescribed Densities.}

\medskip

Assume, as in the previous section, that $F$ is a cone subequation with a finite 
Riesz characteristic $p$. For an arbitrary $F$-subharmonic function $u$, the density of 
$u$ at a point $x_0$ in its domain, is the limit
$$
\eqalign
{
&\T(u, x_0) \ \equiv \ \lim_{r\to0}  {\sup_{B_r(x_0)} u  \over K_p(x-x_0)} \qquad \qquad\   {\rm for} \ \ p\geq 2
\qquad{\rm and}   \cr
& \T(u, x_0) \ \equiv \ \lim_{r\to0}  {\sup_{B_r(x_0)} u  - u(x_0) \over K_p(x-x_0)} \quad {\rm for} \ \ 1\leq p<2,
}
\eqno{(\CC.1)}
$$
which always exists by   [\HLTangI].

The {\bf Dirichlet Problem with Prescribed Densities},  denoted by (DPPD), is the same as the (DPPS)
except that the asymptotic requirement
$$
H \ \approx\ \T_j K_p(x-x_j) \quad {\rm for} \ p\geq2
\qquad{\rm or}\qquad
H \ \sim   \ \T_j K_p(x-x_j) \quad {\rm for} \ 1\leq p<2
\eqno{(\CC.2)}
$$
is replaced by prescribing the density
$$
\T(H, x_j) \ =\ \T_j
\eqno{(\CC.3)}
$$
at each point $x_j$.

\Lemma{\CC.1} {\sl
For any $F$-subharmonic function $u$, the  condition (\CC.2) implies (\CC.3).}
\pf
Assume $x_j=0$.  The notion $u\sim \T K_p$ at the origin is defined for $p\geq2$ by requiring that
$$
\lim_{x\to0} {u(x)\over K_p(x)} \ =\ \T
$$
(see (A.1) and also (4a) in (DPPS)). By Proposition A.3
the condition $u\approx \T K_p$ at 0 implies that $u \sim \T K_p$ at 0, so even when
$p\geq2$, we can assume $u\sim \T K_p$.  Finally, the limits in (\CC.1) equal $\T$ since the limits in 
(A.1) as $x\to 0$ equal $\T$.\qed

\Cor{\CC.2}
{\sl
If existence holds for the (DPPS), then existence holds for the (DPPD). 
Moreover, uniqueness holds for the (DPPD) if and only if for any downward-pointing
singular $F$-harmonic $H$ with density $\T(H,0) = \T$ we have 
$$
H\ \approx\ \T K_p\ \ {\rm at\ } 0 \ \ \ {\rm if \ } 2 \leq p
\qquad{\rm or}\qquad
H\ \sim\ \T K_p\ \ {\rm at\ } 0   \ \ \  {\rm if \ } 1\leq p<2
\eqno{(\CC.4)}
$$
}
\bigskip

\centerline{\bf The Polar Case $(p\geq2)$}

\medskip

Existence for the (DPPS) provided by Theorem \BB.5 implies existence for the (DPPD) because of Corollary \CC.2.

\medskip
\noindent
{\bf THEOREM \CC.3. (Polar Case $p\geq 2$)} \  {\sl Let $F$, $\O$ and $\vf$
be as in Theorem \BB.5.  Then there exits a  
solution $H$ to the (DPPD)
with prescribed positive density $\T$ at any given point $x_0\in \O$.
Moreover,  if $F$ is convex, there exists a solution $H$ with prescribed positive
densities $\T_1, ... , \T_k$ at an arbitrary collection of distinct points $x_1, ... , x_k \in\O$.
 }

\medskip
Uniqueness of the non-linear Green's function implies uniqueness in the general (DPPD).

\Theorem{\CC.4. (Polar Case $p\geq 2$)} 
{\sl 
Suppose that the Riesz kernel $K_p$ is the only solution to the (DPPD) on a ball about the origin with
the same boundary values and asymptotic behavior as $K_p$.  Then
uniqueness holds for the general Dirichlet problem with prescribed densities   (DPPD)  in Theorem \CC.3.}

\pf
Suppose $H$ is a solution to the (DPPD) in Theorem \CC.3. 
By the definition of the (DPPD) we know that $H$ has density $\T_j$ at $x_j$.
It suffices to prove  (\CC.4), i.e., to show that $H\approx \T_j K_p(x-x_j)$
because then the uniqueness part of Theorem \BB.5 applies.
As in [La$_3$, Thm. 3.6] here is how   (\CC.4) can be proved for $H$.
Normalize so that $x_j=0$, $\T_j=1$ and $B\ss\ss\O$ is a ball about the origin.
By the existence part of Theorem \CC.3, we obtain $h\in C(\overline B -\{0\})$, which
is $F$-harmonic on $B-\{0\}$, equal to the constant $K_p$ on $\partial B$, and
satisfies $h \approx H$ at 0.  Now $h \approx H$ at 0 implies that $h$ also has density 1 at 0.
By the hypothesis, this proves $h=K_p$.\qed

\Cor{\CC.5. (Polar Case $p\geq 2$)}
{\sl
Suppose $F$ is an orthogonally invariant subequation of finite Riesz characteristic $p$.
Then both existence and uniqueness hold  for the (DPPD) in Theorem \CC.3.
}
 
\pf
 The classical moving plane argument shows that if $h$ is $F$-harmonic on $B$, constant on $\partial B$ and
 $h(0)=-\infty$, then $h$ is a radial function
(see [La$_3$], [GLN]).  Hence, $h(x) = \T K_p(x)+k$  by  [\HLTangI, Prop. 3.5].\qed

\medskip
\noindent
{\bf Final Note.}  For orthogonally invariant subequations $F$ (and many others as well) it has
been shown ([\HLTangI], [\HLTangII]) that for any $F$-{\bf sub}harmonic function $u$
on $\O\ss\rn$ and $c>0$,
$$
{\sl the\  set\ \ } E_c(u)\ \equiv\ \{x\in \O : \Theta(u,x) \geq c\}\ \ {\sl is \ discrete.}
$$  
The results here show that any
finite set $E\ss \O$ can occur as $E_c(H)$ for an $F$-harmonic function $H$ on $\O$.

\vskip.3in


\centerline{\headfont
\DD.  Comparison.
}
\medskip

In this section we prove the Uniqueness Theorem \AAA.3 for the (DPPS) 
by establishing a comparison theorem in the setting of prescribed singularities.  
The key idea is contained in a local result. First we consider the polar case.
 
\Lemma{\DD.1. (The Polar Case)}  
{\sl
Suppose $v$ is $F$-subhamronic  and $w$ is $\ft$-subharmonic in a deleted neighborhood of 
a point $x_0\in\rn$, and that $h$ is a downward-pointing singular $F$-harmonic at $x_0$. Assume
that for some constants $c$ and $k$
$$
v\ \leq \ h+c \and w\ \leq \ -h +k
\qquad{\rm near}\ \ x_0.
\eqno{(\DD.1)}
$$
Then, with $u\equiv v+w$ extended to $x_0$ by setting
$$
u(x_0) \ \equiv \ \overline{\lim_{x\to x_0}} (v(x) + w(x))
\eqno{(\DD.2)}
$$
we have that}
$$
u\ \ {\sl is\ subaffine\ on\ a\  neighborhood\ of\ \ } x_0.
\eqno{(\DD.3)}
$$
\pf We recall the notion of subaffine functions [\HLDD].  These are the functions $u\in\USC(\O)$ such that
 for any affine function $a$ and  $K\ss\ss \O$,\ \ 
  $u\leq a$ on $\partial K \ \ \Rightarrow\ \  u\leq a$ on $K$.
It turns
out that this is a local condition on the function, namely that it satisfy the subequation $\cpt$ dual to
$\cp$ (see [\HLDD]).
Comparison follows from (\DD.3) since subaffine functions clearly satisfy the maximum principle.

Consider the sum $u=v+w$ which defines an upper semi-continuous $[-\infty, \infty)$-valued function
in a deleted neighborhood of $x_0$.  This function has the following two properties:
$$
u\ \ {\rm is\ bounded\ above \ across\ \ } x_0,
\eqno{(\DD.4)}
$$
$$
u\ \ {\rm is\ subaffine \ on\ a \  deleted\ neighborhood \ of\ \  } x_0.
\eqno{(\DD.5)} 
$$
Obviously (\DD.1) implies (\DD.4).  For (\DD.5) recall (see [\HLDD] or [\HLAE, Thm. 6.2]) that for any 
constant coefficient, pure second-order subequation
$F$, the sum of an $F$-subharmonic function and an $\ft$-subharmonic function is $\cpt$-subharmonic.
The condition (\DD.4) implies that the extension of $u$ defined by (\DD.2) is upper semi-continuous,
with values in $[-\infty, \infty)$, in a neighborhood of $x_0$.  It remains to show that $u$ is $\cpt$-subharmonic.
(Note that the function $ -|x-x_0|$ satisfies (\DD.4) and (\DD.5) and, even though it is continuous across $x_0$,
it is not $\cpt$-subharmonic across $x_0$.  Said differently, Lemma \DD.1 is not simply a removable 
singularity theorem for $\cpt$.)

To prove that $u$ is subaffine across $x_0$ we approximate $u$ by
$$
u_\e \ =\ u + \e v \ =\ (1+\e) v + w.
\eqno{(\DD.6)}
$$
Since $u$ is bounded above and $v(x_0) = -\infty$, if we define $u_\e(x_0) = -\infty$,
then $u_\e$ is upper semi-continuous on a neighborhood of $x_0$.  
Note that $u_\e$ has no test functions at $x_0$, so to prove
that $u_\e$ is subaffine on a  neighborhood $V$ of $x_0$, we need only prove that 
$u_\e$ is  subaffine on $V-\{x_0\}$. However, $u_\e$ is the sum of the $F$-subharmonic
function $(1+\e)v$  and the $\ft$-subharmonic function $w$. This implies, as above,  that
$u_\e$ is $\cpt$-subharmonic on $V-\{x_0\}$ as desired.
(Here we have used that $F$ is a cone.)

Now on a neighborhood of $x_0$ the subaffine functions $u_\e$ increase pointwise to $u$ (since $v<0$).
By the ``families bounded above'' property (see [\HLDD]) and the fact that $u$ is upper semi-continuous
with values in $[-\infty, \infty)$, this proves that $u$ is $\cpt$-subharmonic.\qed
\medskip

Comparison can be stated as follows. 

\Theorem{\DD.2. (Comparison in the Polar Case)} {\sl
Suppose $\O$ is a domain and $h_1,...,h_k$ are downward-pointing singular $F$-harmonics at
$x_1,...,x_k$ respectively.  Given $v,w\in \USC(\ob-\{x_1, ... , x_k\})$
with $v$ $F$-subharmonic and $w$ $\ft$-subharmonic on 
$\O-\{x_1, ... , x_k\}$, suppose that near each $x_j$, $j=1,...,k$ we have
$$
v\ \leq \ h_j + c_j \and w\ \leq \ -h_j + k_j \qquad {\sl for\ some\  constants\ \ } c_j\ {\sl and\ \ } k_j.
\eqno{(\DD.7)}
$$
Then comparison holds on $\O$, that is,
$$
{\sl If\ \ } v+w\ \leq \ 0 \ \ {\sl on}\ \ \bo, \ \ {\sl then\ \ } v+w\ \leq\ 0  \ \ {\sl on}\ \ \O-\{x_1, ... , x_k\}.
\eqno{(\DD.8)}
$$
}
\pf
By Lemma \DD.1 the function $u \equiv v+w$, defined on $\ob-\{x_1, ... , x_k\}$, 
extends to an upper semi-continuous $[-\infty, \infty)$-valued function $u$ on $\ob$
which is subaffine on $\O$.  Hence, $\sup_{\ob} u \leq \sup_{\bo} u$
by the maximum principle.\qed

\medskip

\noindent
{\bf Proof of the Uniqueness Theorem \AAA.3 in the Polar Case.}
It suffices to prove that 
$$
{\rm If}\ \ H\ \ {\rm is\ a \ solution\ to\  the\ (DPPS), \ 
then\ }\ H(x) \ =\ \sup_{v\in \cf} v(x).
\eqno{(\DD.9)}
$$
where $\cf$ is the family defined in Theorem \AAA.4.
We can apply the Comparison Theorem \DD.2 to $v\in \cf$ and $w\equiv -H$.  On $\partial \O$,
we have $v\leq \vf$ and $w  = -\vf$, and so $v+w\leq0$ (or $v\leq H$) on $\ob$ by (\DD.8). Since
$H\in \cf$, this proves (\DD.9).\qed

\medskip

Now we turn to the finite case.

\Lemma{\DD.3. (The Finite Case)}
{\sl
Suppose $v$ is $F$-subharmonic and $w$ is $\ft$-subharmonic on a deleted neighborhood of $x_0$,
and that $h$ is a downward-pointing singular $F$-harmonic at $x_0$.  Assume that 
$$
\liminf_{x\to x_0} {v(x)-v(x_0)  \over h(x) - h(x_0)}\ \geq\ 1
\and
\liminf_{x\to x_0} {w(x)-w(x_0)  \over h(x) - h(x_0)}\ \geq\ -1
\eqno{(\DD.10)}
$$
Then with $u\equiv v+w$ extended to $x_0$ by setting
$$
u(x_0) \ \equiv  \ \limsup_{x\to x_0} u(x)
\eqno{(\DD.11)}
$$
we have that 
$$
u\ \ {\sl is\ subaffine\ on\ a \ neighborhood\ of\ \ } x_0.
\eqno{(\DD.12)}
$$
}

\pf
To prove that $u$ is subaffine, we approximate $\overline u  \equiv u - v(x_0)+w(x_0)$ by
$$
{\overline u}_\e \ =\ \overline u + \e(v(x) -v(x_0)) \ =\ (1+\e) (v(x)-v(x_0)) + (w(x)-w(x_0))
\eqno{(\DD.13)}
$$
and prove that  $\overline u$ is subaffine.  Note that as in Lemma \DD.1, 
${\overline u}_\e$ is subaffine on a deleted neighborhood of $x_0$ since $F$ is a cone.
To show that ${\overline u}_\e$ is subaffine on a neighborhood of $x_0$, we need only show that at
$x_0$, ${\overline u}_\e$ has no test functions.

The hypothesis (\DD.10) on $v$ implies that for $1<\a<1+\e$, there exists a neighborhood
of $x_0$ with $v(x) -v(x_0) \geq {\a\over 1+\e} (h(x)-h(x_0))$.

The hypothesis (\DD.10) on $w$ implies that for $1<\b<\a$,  we have 
$
w(x)-w(x_0) \geq (\a-\b) (h(x)-h(x_0))
$
near $x_0$.  Hence,  ${\overline u}_\e \geq  (\a-\b) (h(x)-h(x_0))$, which proves that 
${\overline u}_\e$ has no test functions at $x_0$, since $h$ has no test functions at $x_0$
(note that ${\overline u}_\e (x_0)=0$.)
This proves that each ${\overline u}_\e$ is subaffine in a  neighborhood of $x_0$.

Finally the fact that $v(x)-v(x_0)$ is bounded below by a positive multiple of $h(x)-h(x_0)$
implies that ${\overline u}_\e$ is decreasing pointwise as $\e\to0$ in a neighborhood of $x_0$.
However, outside $x_0$, $\overline u$ is the sum of the $F$-subharmonic function $v(x)-v(x_0)$
and the $\ft$-subharmonic function $w(x)-w(x_0)$, which implies that $\overline u$ is subaffine.\qed

\Theorem{\DD.4. (Comparison in the Finite Case)}
{\sl
Suppose that $\O$ is a domain and $h_1,...,h_k$ are downward-pointing singular $F$-harmonics at
$x_1,...,x_k \in\O$ respectively.  
Given  $v,w \in \USC(\ob - \{x_1,...,x_k\})$   
with $v$ $F$-subharmonic and $w$ $\ft$-subharmonic on 
$\O-\{x_1, ... , x_k\}$, suppose that  for $j=1,...,k$ we have
$$
\liminf_{x\to x_j} {v(x)-v(x_j)  \over h(x) - h(x_j)} \ \geq \ 1
\and
\liminf_{x\to x_j} {w(x)-w(x_j)  \over h(x) - h(x_j)} \ \geq \ -1.
\eqno{(\DD.14)}
$$
Then (with $u\equiv v+w$ defined at $x_j$ by $(v+w)(x_j) \equiv \limsup_{x\to x_j} (v+w)(x)$),
$$
{\sl If\ \ } v+w\ \leq \ 0 \ \ {\sl on}\ \ \bo, \ \ {\sl then\ \ } v+w\ \leq\ 0  \ \ {\sl on}\ \ \ob.
\eqno{(\DD.15)}
$$
}
\pf
This follows from Lemma \DD.3 as in the proof of Theorem \DD.2.\qed

\Cor{\DD.5. (Uniqueness in the Finite Case)}
{\sl
Let $\cf$ denote the family of $v\in\USC(\ob)$ satisfying: $v$ is $F$-subharmonic on $\O$, $v\leq \vf$ on $\bo$,
and for $j=i,...,k$
$$
\liminf_{x\to x_j} {v(x)-v(x_j)  \over h(x) - h(x_j)} \ \geq \ 1.
$$
If $H$ is a solution to the (DPPS), then 
$$
H(x) \ =\ \sup_{v\in \cf} v(x).
\eqno{(\DD.16)}
$$
}
\pf
Note that if $v\in\cf$ and $H$ is a solution, then comparison implies that $v\leq H$.
Since $H\in \cf$, this proves (\DD.16).\qed

\vskip .3in


\def\vfl{\underline{\vf}}
\def\hh{h}

\centerline{\headfont \AA.\   A Basic Construction and the Proof of Existence in the Polar Case.}
\bigskip

In this section we describe the basic construction of the solution to the (DPPS) in some generality,
and then complete the existence proof in the polar case.
The finite case  is finished in Section \EE. Our starting point is the standard (DP) on the domain
with a neighborhood of the singular points removed.

\medskip
\centerline
{\bf The Dirichlet Problem on the Perforated Domain.}
\medskip

We  fix $r_0>0$ so that the closed balls $\overline B_{r_0}(x_j) = \{|x-x_j|\leq r_0\},\  j=1,...,k$, are mutually disjoint and contained in $\O$.  Let 
$$
D_r \ \equiv \  \bigcup_{j=1}^k  \overline B_r(x_j) \qquad{\rm for}\ 0<r\leq r_0.
$$
Consider the perforated  domain
$$
\O_r \ \equiv \  \O - D_r
$$
whose oriented boundary   is the sum of the {\sl outer boundary} $\bo$ and the {\sl  inner boundary} $-\partial D_r$.

 Consider the Dirichlet Problem for a subequation $F \ss\Symn$
on  $\O_r$
with given boundary functions $\vfl \in C(\partial D_r)$ 
and $\vf \in C(\bo)$. 
As discussed in the proof of Lemma \DD.1 comparison and hence 
uniqueness holds on any domain. For existence,
consider the {\bf Perron family} $\cf$ consisting of all $u\in\USC( \overline{ \O_r})$
such that
$$
u\bigr|_{\O_r} \in F(\O_r),  \qquad
u\bigr|_{\partial D_r} \leq \vfl, \and u\bigr|_{\bo} \leq \vf.
$$
along with its {\bf Perron function}
$$
H_r(x) \ \equiv\ \sup_{u \in \cf} u(x).
$$

By [\HLDDR, \S 12] the Perron function will solve the Dirichlet problem for the given boundary values 
provided that,  at each point of  $\partial \O_r = \partial \O - \partial D_r$,  one can construct barriers
 as in Propositions $F$ and $\ft$ in  [\HLDDR]  on page 453.
 
 We list our assumptions.
 
 \medskip
 \noindent
 {\bf Assumption (B1):}  The outer boundary $\bo$ is strictly $F$-convex.
 
 \medskip
 
 Note that $\bo$ is strictly $\ft$-convex by (F3) = (B3). Thus  barriers, as in Propositions $F$ and $\ft$
 exist for each point  $x_0\in \partial \O$.  

Note similarly that  the inner boundary $-\partial D_r$ is also strictly $\ft$-convex by (F3) = (B3). 
   This provides a barrier as 
in Proposition $\ft$ at points $x_0\in \partial D_r$.  

To obtain an $F$-barrier  at a point $x_0\in\partial D_r$ we 
assume that the inner boundary function $\vfl$ is of a special nature, namely,

 \medskip
 \noindent
 {\bf Assumption (B2):} $\vfl \equiv \pss\bigr|_{\partial D_r}$ where  $\pss$ is $F$-subharmonic on $\ob$, and either

 \medskip
 {\bf Polar Case:}  $\pss:\ob \to [-\infty, \infty)$ is continuous  and $= -\infty$ precisely
 at the points 
 
  \hskip .9in $x_1, ... , x_k$,  or

 \medskip
 
 {\bf  Finite Case:}  $\pss \in C(\ob)$ with each  $x_j$, $j=1,...,k$ a strict local minimum point

 \hskip .9in      with no test functions.
 
 \medskip


 \medskip
 \noindent
 {\bf Assumption (B4):}  \ \ $\pss\bigr|_{\bo} \ \leq\ \vf$.
 
 \medskip

Since $\pss$ and $\vf$ satisfy (B2) and (B4), we can  use the $F$-subharmonic function $\pss$ to construct  barriers at points $x_0\in\partial D_r$
as in Proposition $F$.
This is done by setting   $\underline u(x) \equiv \pss(x) -\d +\e|x-x_0|^2$ for $\e>0$ sufficiently small. 

This establishes the following existence result.

\Theorem{\AA.1}
{\sl
Assume that $F$ is a cone subequation which satisfies Condition (F3) = (B3).
Suppose $\vf \in C(\partial \O)$.
Suppose that $\bo$ is strictly $F$-convex  (B1),
and that  $\pss$ satisfies (B2) and (B4). 
Then the Perron function $H_r$ solves the 
Dirichlet Problem:

\medskip

(a) \ \ \ $H_r \in C(\ob_r)$

\medskip

(b) \ \ \ $H_r$  is $F$-harmonic on $\O_r$

\medskip

(c) \ \ \ $H_r \bigr|_{\partial D_r} \ =\ \pss\bigr|_{\partial D_r} 
 \and 
 H_r \bigr|_{\bo} \ =\ \vf$.
}

\vskip.3in
\centerline{\bf The Candidate for the Solution to the (DPPS)}

\medskip

We continue with the same notation and hypotheses.

The proposed solution  to the (DPPS)  is constructed as a pointwise increasing limit 
$$
H(x) \ \equiv\ \lim_{r\to0} \overline H_r(x) \qquad {\rm on}\ \ \ob 
\eqno{(\AA.1)}
$$
of functions $\overline H_r$  which extend the functions $H_r$ to $\ob$.

\Lemma{\AA.2}  {\sl
The  continuous functions
$$
\overline{H}_r \ \equiv \ 
\cases
{
H_r \quad {\rm on}\ \ \ob_r  \cr
\ \pss  \quad\ \   {\rm on}\ \  D_r
}
\eqno{(\AA.2)}
$$
are $F$-subharmonic on $\O$ and pointwise increasing on $\ob$ as $r\to0$.}
\pf
The function $\pss$ is in the Perron family on $\O_r$ since $\pss\bigr|_{\partial D_r} =\vfl$ 
and $\pss\bigr|_{\bo} \leq \vf$. This proves that 
$$
\pss\leq H_r \quad {\rm on}\ \  \ob_r.
\eqno{(\AA.3)}
$$
Consequently,
$$
\overline H_r(x) \ \equiv\ \cases
{
\max\{ \pss(x), H_r(x)\} \qquad {\rm on}\ \  \ob_r  \cr
\ \ \qquad \pss(x) \qquad\quad\qquad\  {\rm on}\ \ D_r.
}
\eqno{(\AA.4)}
$$
To see  that
$$
\overline H_r \ \ \ {\rm is\ } F\,{\rm subharmonic\ on\ \ }\O,
\eqno{(\AA.5)}
$$
 note that $\overline H_r^\e = \max\{\pss+\e, H_r\}$ defines an $F$-subharmonic 
function on $\O$ since $\pss+\e >H_r$ on $\partial D_r$,  and  then note that $\overline H_r^\e$ decreases pointwise 
as $\e\to0$ to $\overline H_r$ on $\ob$.  

Finally we show that for all $\rho$ with $0<\rho < r\leq r_0$, one has
$$
\overline H_r\ \leq\ \overline H_\rho\qquad{\rm on}\ \ \overline \O.
\eqno{(\AA.6)}
$$
On $D_r$ we have $\overline H_r = \overline H_\rho = \pss$.  
In particular, on $\partial D_\rho$, $\overline H_r = \pss$, while on 
the outer boundary $\bo$, $\overline H_r =\vf$.
Thus $\overline H_r$ is in the Perron family for  $H_\rho$ on this larger domain $\O_\rho$, 
proving that $\overline H_r\leq H_\rho$ on $\ob_\rho$, which establishes (\AA.6).\qed

Let $H^{\rm DP}$ denote the solution to the standard Dirichlet Problem (DP) on $\O$.
That is, $H^{\rm DP}\in C(\O)$, $H^{\rm DP}$ is $F$-harmonic on $\O$, and 
$H^{\rm DP}\bigr|_{\bo}=\vf$.  This function exists since the outer boundary is assumed
to be strictly $F$-convex, and by (F3) = (B3) it is also strictly $\ft$-convex.

\Prop{\AA.3} {\sl
The proposed solution $H$ to the (DPPS) defined by (\AA.1) satisfies:

\medskip
(1) \ (a)\ \  $H^* \in \USC(\ob)$ \and  (b) \ \ $-H \ =\ (-H)^* \in \USC(\ob-\{x_1, ... , x_k\})$,

\medskip
(2) \ (a)\   $H^*$ is $F$-subharmonic on $\O$, 
\ \ \   
(b)  \ $-H $ is $\ft$-subharmonic on $\O-\{x_1, ... , x_k\}$

\medskip
(3) \  $H^*\bigr|_{\bo} \ =\ H\bigr|_{\bo} \ =\ \vf$

\medskip
(4) \ \  $\pss \ \leq\ H \ \leq\ H^{\rm DP} \ \ \ $ on $\ob$.
}
\pf
Assertion (1) is immediate while (2a) follows from the 
``families bounded above property'' and (2b) by the  ``decreasing limit property''.
 To prove (4) first note that (\AA.3) implies that $\pss\leq H$.
 Now the function
  $H^{\rm DP}$ is the Perron function for the standard Dirichlet Problem on $\O$ with boundary values
$\vf$.  Therefore, $H^{\rm DP}$ is larger than $\pss$ by assumption (B4).  Thus $H_r \leq H^{\rm DP}$ on $\bo_r$,
and this remains true on $\ob_r$.  Therefore, $H\leq H^{\rm DP}$, which  implies that 
$$
H^*\ \leq \ H^{\rm DP}.
\eqno{(\AA.7)}
$$
In turn this implies (3).  \qed

\bigskip
\centerline
{\bf
Existence in the Polar Case -- The Proof of Theorem \AAA.4
}
\medskip

First note that in the polar case the notion of asymptotic equivalence is preserved
by subtracting a constant. This implies that the hypothesis (B4) in Theorem \AA.1
and Proposition \AA.3 can always be satisfied.

Now (4) in Proposition \AA.3 yields the left hand inequality in Part (4a) of the 
(DPPS) stated in the introduction.  For the right hand inequality in (4a) we prove
the following.

\Lemma {\AA.4} {\sl
Near each $x_j$ we have $H\leq \hh_j + c_j$ for some constant $c_j$.
}
\pf
For $0< \rho\leq r$ we have that 
$$
H_\rho -\hh_j \quad {\rm is}\ \ \cpt\ {\rm subharmonic\ on\ \ } A_{\rho,r} \ =\  B_r  -  \overline B_\rho.
$$
On the inner boundary  $\partial B_\rho$ we have, since $\pss \approx \hh_j$ implies $\pss -\hh_j\leq C$, that 
$$
 H_\rho - \hh_j \ =\   \pss -\hh_j \ \leq \  C,
$$
and on the outer boundary $\partial B_r$ we have
$$
 H_\rho - \hh_j \ \leq \   U-\hh_j \ \leq\ C(r)
$$
independent of $\rho$.  Hence,
$$
H- \hh_j\ \leq\ \max\{C,C(r)\} \ \equiv\ C'\qquad{\rm on}\ \ B_r(x_j). \qquad\qquad\mathqed
 $$
  
 \medskip
 
Now we apply the Comparison Theorem \DD.2 to $v\equiv H^*$ and $w\equiv - H$.
By (1a) $H^*$ takes values in $[-\infty, \infty)$ on $\O$, while by (4) 
$w$ takes values in  $(-\infty,  \infty)$ except at $x_1,...,x_k$ where $w$ equals $+\infty$.
The inequality 
$$
v \ \equiv \  H^* \ \leq\ \hh_j + c_j    \qquad{\rm near}\ \ x_j
$$
is immediate from Lemma \AA.4.  Observe now that by combining (4) $\pss  \leq H \ (\equiv -w)$
from Proposition \AA.3  with  the inequality $\hh_j - k_j \leq \pss$ near $x_j$, which is part of Hypothesis 
(H) in Theorem \AAA.4,  gives
$$
w\ \leq \ -\hh_j + k_j   \qquad{\rm near}\ \ x_j.
$$
This establishes the hypothesis (\DD.7) in Theorem \DD.2.
Finally, by (3) in Proposition \AA.3, $v=\vf$ and $w=-\vf$ on $\bo$, so that $v+w = 0$ on $\bo$.
We can now apply  Theorem \DD.2 to conclude that 
$v+w \leq 0$ on $\ob - \{x_1,...,x_k\}$, i.e.,  $H^*\leq H$ on $\ob - \{x_1,...,x_k\}$.
This proves that $H^*=H$ on  $\ob - \{x_1,...,x_k\}$.  Since we already have 
$H=H_*$ (because it is an increasing limit of continuous functions), we conclude
that $H$ is continuous on $\ob - \{x_1,...,x_k\}$, and Conditions (1) and (2) are proved.
This completes the proof of Theorem \AAA.4.  \qed

\vfill\eject
\centerline
{\bf
The (DPPS) with Singularities on a Compact Polar Set
}

\medskip
 
The arguments given above adapt  to prove a version of Theorem \AAA.4
with $\{x_1, ... , x_k\}$ replaced by a compact polar set.

\Theorem{ \AA.5} 
{\sl
Let $F$ be a  cone subequation  satisfying  Condition  (F3).  
Let $\O\ss\ss \rn$ be a domain with smooth boundary
$\bo$ which is strictly $F$-convex, and let $\Sigma\ss \O$ be a compact subset.
Suppose there exists a continuous function  $h:\ob\to [-\infty, \infty)$ such that 
$\Sigma = h^{-1}(-\infty)$ and $h$ is $F$-harmonic on $\O-\Sigma$.

  Then for any $\vf\in C(\bo)$ we have the following.
\medskip
\noindent
{\bf Existence.}  There exists $H\in C(\ob - \Sigma)$ such that:
\medskip

(1)  \quad  $H$ is $F$-harmonic on $\O -\Sigma$,
\medskip

(2)  \quad  $H\bigr|_{\bo} = \vf$,
\medskip

(3)  \quad  $H$ is asymptotically equivalent  to $h$,
i.e.,  there exists $c, C\in \bbr$ such that
$$
h(x) + c \ \leq\ H(x)\ \leq\  h(x) +C  \qquad  {\rm on}\ \  \O-\Sigma
$$
\medskip

\noindent
{\bf Uniqueness.}   There is at most one function $H \in C(\ob-\Sigma)$ satisfying (1), (2) and (3).
(Here the $F$-convexity of $\bo$ is not required.)
}

\pf
The uniqueness is proved as in Section  \DD.
For existence we choose a sequence of regular values $\{r_j\}$ of $h$ with $r_j \downarrow -\infty$.
Again by adjusting $h$ with an additive constant we can assume that
(B4) is satisfied.  For each $j$ define
$$
D_{j} \ \equiv\ \{x\in \O :  h(x) > r_j\} 
\and
\O_j \ \equiv \ \O - \overline D_j.
$$
Let  $H_j$ be the solution to the Dirichlet Problem on $\O_j$ for
$F$-harmonic functions with boundary values 
$$
H_j\bigr|_{\bo}\ =\ \vf \ 
\and 
H_j\bigr|_{\partial D_j}\ =\ h\bigr|_{\partial D_j}
$$
on the outer and inner boundaries respectively.
As before this solution exists due to the assumption of $F$-convexity on
$\bo$, the hypothesis (F3) = (B3),  the fact that $h$ is in the Perron family by (B4),
and the $F$-harmonicity of $h$ outside $\Sigma$.

Now define $\overline{H}_{j}(x)$ as in Lemma \AA.2 and note that the assertions
of Lemma \AA.3    hold by exactly the same arguments.

It follows that the increasing limit 
$$
H\ \equiv\ \lim_{j\to\infty} \overline H_j
$$
has the properties that $H^*$ is $F$-subharmonic  and $-H$ is $\ft$-subharmonic 
on $\O-\Sigma$.  Using the solution $H^{\rm DP}$ to the Dirichlet Problem on $\O$
with boundary values $\vf$ the same arguments as above show that 
$$
H^*\bigr|_{\bo}=\vf.
\eqno{(\AA.8)}
$$

We now prove (3).
The left-hand inequality in (3) holds since $h$ is in the Perron family 
for $H_j$ on $\O_j$ for all $j$. On the other hand, for $C>0$
sufficiently large we will have $h+C >\vf$ on $\bo$. Since $h+C$
is $F$-harmonic on $\O-\Sigma$  and greater than $H_j$ on $\partial \O_j$,
we have $h+C > H_j$ on $\O_j$ by comparison. This establishes the right-hand inequality in (3).

As noted above we have
$$
H^*\in F(\O-\Sigma)\and
-H \in \ft(\O-\Sigma).
$$
As before this implies that $H^*-H\in\cpt_p(\O-\Sigma)$.
Now Condition (3) implies Condition (3) for $H^*$ and therefore
$H^*-H\leq  C-c$ on $\O-\Sigma$.  We now apply the removable singularity
argument in [\HLRS] using the polar  function $h$ (as before) to conclude
that $H^*-H$ is $F$-subharmonic on $\O$. By (\AA.8)  $H^*-H$ is $\leq 0$ on $\bo$.
Hence, we have $H^*-H\leq0$ on $\O$. 
We have proved that $H^*=H$ on $\O$.  This proves $H\in C(\ob-\Sigma)$ and
condition (1), and we are done. \qed

\vskip .3in


\centerline{\headfont \EE.   Existence in the Finite Case.}
\bigskip

We now take up the proof of existence in the finite case.
It proceeds in two stages.  The first (Theorem \EE.1)  is a construction which
provides a family of $F$ harmonics which only satisfy a weakened form of asymptotic equivalence
in a range.  The second stage shows that one member of the family has the desired asymptotic 
singularity at the given point.

Throughout this section we assume that $F$ is a cone subequation with the property 
that  all boundaries are  strictly $\ft$-convex, and that $\O$ is a domain with smooth strictly
$F$-convex boundary. We limit the discussion to the case of a single singular point $x_0\in\O$.
(However, some of the arguments extend to multiple singular points as in the polar case.)
 The following assumption replaces the hypothesis (H) in Theorem \AAA.4.

\medskip
\noindent
{\bf Hypothesis (H1) (Finite Case).}  We are given a function $\ps\in C(\ob)$ which is $F$-harmonic
on $\O-\{x_0\}$ and has a downward-pointing singularity at $x_0$
with $h(x_0) < \inf_{\bo} h$.
\medskip

Applying  the maximum principle to $-h$  on  $\O-B_r(x_0)$, for $r$ small, proves  that:
$$
{\rm the\  point}\ \ x_0 \ \ {\rm   is \ a \  strict\  global \  minimum \ for \ } \  h\ \ {\rm on\ }\  \ob.
\eqno{(\EE.1a)}
$$
Moreover, for convenience, 
we assume by rescaling  that 
$$
\ps(x_0) \ =\ 0, \qquad  \sup_{\bo} \ps\ =\ 1
\and 
h(x) \ >\ 0\ \ {\rm for}\ \ x\neq x_0\ \ {\rm in}\ \ob.
\eqno{(\EE.1b)}
$$

\bigskip
\centerline
{\bf The Construction}
\medskip

It is similar to the construction in the polar case.   However, for each $t\geq0$  we construct a function
$H^t$ which is a candidate for the solution with $H^t \sim t\ps$ at $x_0$.
When $t=0$, the construction will yield the solution to the standard Dirichlet Problem on $\O$, namely:
$$
\Hz \in C(\ob), \ \ \Hz \ {\rm is \ } F{\,\rm harmonic\  on\ } \O, \ \ {\rm and}\ \ \Hz\bigr|_{\bo} \ =\ \vf.
\eqno{(\EE.2)}
$$

The construction of $H^t$ is based on using the function
$$
\hlb_t (x) \ \equiv\ \Hz (x) + t(\ps(x)-1).
\eqno{(\EE.3)}
$$
to prescribe the boundary values on the inner boundary $\partial B_r(x_0)$.
An upper bound for $H^t$ is provided by the function
$$
\hub_t (x)\ \equiv\  \l \ps(x) + \Hz (x_0) + t(\ps(x)-1)
\eqno{(\EE.4)}
$$
for $\l$ sufficiently large.  Note that $\hub_t$ is $F$-harmonic on $\O-\{x_0\}$
since $h$ is $F$-harmonic there.  We need the following  ``monotonicity'' hypothesis
to ensure that $\hlb_t$ is $F$-subharmonic.

\medskip
\noindent
{\bf Hypothesis (H2).}  The function $\Hz + t\ps$ is $F$-subharmonic on $\O$, and therefore, so is each  $\hlb_t$.
(Of course (H2) is satisfied if $F$ is a convex subequation.)
\medskip

Let $H^t_r$ denote the solution, given by Theorem \AA.1,
 to the Dirichlet Problem on $\O_r = \O -\overline{B_r(x_0)}$ with 
boundary values $\vf$ on $\bo$ and boundary values $\hlb_t$ on $\partial B_r(x_0)$.
As in Section \AA, we will show that these functions  $H^t_r$ are pointwise increasing  as $r\to0$.
We define $H^t$ by
$$
H^t(x)  \ \equiv\ \lim_{r\to0} H^t_r(x) \ \ {\rm on}\ \ob-\{x_0\} 
\and
H^t(x_0)  \ \equiv\ \Hz(x_0) -t.
\eqno{(\EE.5)}
$$
This is our candidate for a ``solution'' to the (DPPS) with singularity prescribed by
$t\ps$ at $x_0$.  Next we show that $H^t$ satisfies all the required conditions with the
exception of the asymptotic equivalence $H^t \sim t\ps$ at $x_0$.

\Theorem{\EE.1} {\sl For each $t\geq0$,
\smallskip
(1) \ \  $H^t \in C(\ob)$,

\smallskip
(2) \ \  $H^t$ is $F$-harmonic on $\O-\{x_0\}$,

\smallskip
(3) \ \  $H^t \bigr|_{\bo} =\vf$.
\smallskip
\noindent
In addition,
\smallskip
(4a) \ \  $\hlb_t \leq H^t $ with equality  at $x_0$,

\smallskip
(4b) \ \  $H^t \leq \hub_t$ (for $\l$ large) with equality  at $x_0$.
}

\pf
First note that $\hlb_t \in C(\ob)$ has the properties
$$
\hlb_t \ \ {\rm is\ \ } F{\,\rm subharmonic\ on\ } \O
\and
 \hlb_t \leq \vf
\eqno{(\EE.6)}
$$
because of (H2) and  $\ps \bigr|_{\bo} \leq 1$ respectively.
(We also have $\hlb_t(x_0) = \Hz(x_0)-t$ 
since  $\ps(x_0)=0$.)  Therefore $\hlb_t$ is in the Perron family
for $H^t_r$, which proves that
$$
\hlb_t \ \leq\ H^t_r \quad {\rm on}\ \ \O_r.
\eqno{(\EE.7)}
$$
Define
$$
{\overline H}^t_r \ \equiv\ 
\cases
{
\max\{\hlb_t, H^t_r\} = H^t_r \quad{\rm on}\ \ \ob_r  \cr
\qquad\qquad\qquad\quad \hlb_t   \quad{\rm on}\ \ \overline{B}_r.
}
\eqno{(\EE.8)}
$$
Then exactly as in Section \AA \  one proves that
$$
{\overline H}^t_r \ \ {\rm is}\ F{\,\rm subharmonic\ on\ }\O\ {\rm and \ increasing\ in\ } r \ {\rm as}\ r\to0.
\eqno{(\EE.9)}
$$
Define
$$
H^t (x) \ \equiv\ \lim_{r\to0}  {\overline H}^t_r(x) \qquad{\rm for\ all\ \ } x\in \ob.
\eqno{(\EE.10)}
$$
Now Property (3) is immediate since $H^t_r\bigr|_{\bo} = \vf$ independent of $r$.
Next note that  $\hlb_t \leq  {\overline H}^t_r \leq H^t$ by (\EE.7).
Together with the equality $\hlb_t(x_0) =  {\overline H}^t_r(x_0) = H^t(x_0)-t$ at $x_0$, this proves (4a).

We begin the proof of (1) and (2).  The family $\{  {\overline H}^t_r\}_{r>0}$ of functions on $\ob$
is bounded.  The upper bound ${\overline H}^t_r \leq \Hz$ implies that 
$H^t \leq (H^t)^* \leq \Hz$, and hence (3) can be strengthened to 
\smallskip

(3)$'$ \ \ $H^t \bigr|_{\bo}\ =\ (H^t)^* \bigr|_{\bo}\ =\ \vf$.
\smallskip
\noindent
Also, since $H^t$ is an increasing limit of continuous functions, it is lower semi-continuous.
Thus,
$$
(H^t)^* \in \USC(\ob) 
\and
H^t = (H^t)_* \in {\rm LSC}(\ob)
\eqno{(\EE.11)}
$$
are finite-valued.  Moreover, by the ``families bounded above'' property and the ``decreasing limit'' property,
we have
$$
(a)\ \ (H^t)^* \ {\rm is}\ F{\,\rm subharmonic\ on\ }\O,
\quad
(b)\ \ -H^t \ {\rm is}\ \ft{\,\rm subharmonic\ on\ }\O - \{x_0\}.
\eqno{(\EE.12)}
$$
Hence $u\equiv (H^t)^*  - H^t \in\USC(\ob)$ satisfies $u\geq0$ on $\ob$, and  $u\bigr|_{\bo}=0$ by (3)$'$.
As discussed in the proof of Lemma \DD.1, $u$ is subaffine on $\O-\{x_0\}$.  Therefore, by the Maximum Principle
for subaffine functions on the domain $\O-\{x_0\}$ we have
$$
 (H^t)^* (x)  -  H^t(x) \ \leq\  (H^t)^* (x_0)  -  H^t(x_0)\ =\ u(x_0) \quad {\rm on}\ \ \ob.
\eqno{(\EE.13)}
$$

It remains to show that $u(x) \equiv (H^t)^*(x)-H^t(x)$ equals zero at $x_0$, i.e., 
$$
(H^t)^* (x_0)  \ =\ \Hz(x_0) -t.
\eqno{(\EE.14)}
$$
Once this is established, we will have $(H^t)^*=H^t$, which implies both (1) and (2).
Next note that  (\EE.14) follows immediately from the upper bound (4b) since $\hub_t(x_0)= \Hz(x_0)-t$
and $\hub_t\in C(\ob)$.

Thus it remains to prove (4b).
As noted, equality in (4b) holds at $x_0$.  Hence it suffices to prove that $H^t\leq \hub_t$ on $\ob -\{x_0\}$, or that
$$
H^t_r \ \leq \ \hub_t\qquad {\rm on}\ \ \O_r\ \ {\rm for\ small\ \ } r.
\eqno{(\EE.15)}
$$
Since $\hub_t$ is $F$-harmonic on $\O_r$ by (H1), it suffices, by comparison, to show that
$$
H^t_r \ \leq\ \hub_t  \qquad{\rm on}\ \ \partial \O_r.
\eqno{(\EE.16)}
$$
That is, we must show that
$$
\vf \ \leq\ \hub_t
\qquad{\rm on \ the\ outer\ boundary}\ \ \partial \O,\ \ {\rm and}
\eqno{(\EE.16a)}
$$
$$
\hlb_t  \ \leq\ \hub_t
\qquad{\rm on \ the\ inner\ boundary}\ \ \partial B_r.
\eqno{(\EE.16b)}
$$
Since $\inf_{\bo} h >0$ and ${\rm osc}_{\bo} (h) = 1 - \inf_{\bo} h$,
it is straightforward to see that 
$$
\l \ \geq \ {\sup \vf - \Hz(x_0) + t\,{\rm osc}_{\bo} (\ps)  \over \inf_{\bo} \ps}
\quad\iff\quad
\sup \vf \ \leq \ \inf_{\bo} \hub_t
\eqno{(\EE.17)}
$$
which implies (\EE.16a).
We might as well take $\l\equiv\l(t)$ to be the affine function of $t$ defined by equality in (\EE.17).

For (\EE.16b) first note that $\hlb_t \leq H^0$ implies $\hlb_t\bigr|_{\bo} \leq\vf$,
and hence by (\EE.16a), that $\hlb_t\bigr|_{\bo} \leq \hub_t\bigr|_{\bo}$.
We also have $\hlb_t(x_0) = \hub_t(x_0)$ since $h(x_0)=0$.
Since $u\equiv \hlb_t - \hub_t  \in \cpt(\O-\{x_0\})$, we can apply the maximum principle
to $u$ on $\O-\{x_0\}$ to conclude that $u\leq 0$ on $\ob$.
In particular, $u\leq0$ on $\partial B_r(x_0)$, which is (\EE.16b).
This  completes the proof of (4b) and, therefore, of Theorem \EE.1. \qed

\bigskip
\centerline
{\bf Prescribing the Density at a Point}
\medskip

We make two additional assumptions on $\ps$.

\medskip
\noindent
{\bf Hypothesis (H3).}  For each $H$ which is $F$-subharmonic near $x_0\in \O$ and $F$-harmonic
on a deleted neighborhood of $x_0$, 
$$
\lim_{x\to x_0}  {H(x) - H(x_0) \over \ps(x) }  \ \equiv \ \T \qquad {\rm exists \ and\ \ } \T\geq0,
\eqno{(\EE.18)}
$$
that is, $H\sim \T\ps$ at $x_0$ for some $\T\geq 0$.
\medskip

There are many examples where this is true.  They will be discussed later in this section.

\Def{\EE.2}  Under the hypothesis (H3) the {\bf $\ps$-density of $H$ at $x_0$}, denoted $\T_{x_0}(H)$
is defined to be the limit in (\EE.18).
\medskip

The second additional assumption on $\ps$ is that $F$-harmonics have vanishing densities.

\medskip
\noindent
{\bf Hypothesis (H4).}  If $H$ is $F$-harmonic in a neighborhood of $x_0$, then 
$\T_{x_0}(H) =0$, i.e., 
$$
\lim_{x\to x_0} {H(x)-H(x_0) \over \ps(x)} \ =\ 0.
$$.
\medskip

\Def{\EE.3}   The {\bf  Dirichlet Problem with a Prescribed density}, abbreviated (DPPD), is said to be 
{\bf uniquely solvable for $F$} if for all $\vf\in C(\bo)$ and $\Theta \geq0$,
there exists a unique $H\in  C(\ob)$ which is $F$-harmonic on $\O-\{x_0\}$ and satisfies
$$
H\bigr|_{\bo} \ =\ \vf\and \Theta_{x_0}(H)\ =\ \Theta \quad({\rm i.e.\ \ } H\sim \Theta h).
$$

As in Theorem \EE.1, we assume (H1). We shall assume that $F$ is convex, so (H2) is unnecessary. 
 In addition we assume (H3) and (H4).

\Theorem{\EE.4. (Existence for a Prescribed Density)}
{\sl
Assume that $F$ is a convex cone subequation.
Then for each $\T\geq0$ and $\vf\in C(\bo)$ the (DPPD) is uniquely solvable.
}
\medskip

For the terminology  used in the next corollary see the discussion following Corollary \AAA.5.

\Cor{\EE.5.\ (The Nonlinear Green's Function)}
{\sl
There exists a unique nonlinear Green's function $G(x) = G_\O(x; x_0,h)$ for the subequation
$F$ on the domain $\O$ with asymptotic singularity determined by the generalized fundamental
solution $h$.
Moreover, $x_0$ is a strict global minimum point for $G$ on $\O$.
}
\pf By definition $G$ is the unique solution on $\O$ with $h$-density 1 at $x_0$ and boundary values
$\vf\equiv 0$.
Since $G$ has no test functions at $x_0$, $G$ is $F$-subharmonic on $\ob$.
Hence, by the maximum principle, $G\leq \sup_{\bo} G = 0$ on $\ob$.
Now $-G$ is $\ft$-subharmonic on $\O-\{x_0\}$.
Therefore, by the maximum principle applied to $-G$ on $\O-\{x_0\}$,
we first get $G(x_0)<0$ (because $G(x_0)=0$ implies $-G\leq 0$, 
and so $G\equiv 0$ on $\ob$ contradicting $\T_{x_0}(G)=1$).
Then exactly as in the proof of (\EE.1a) we get that $x_0$ is a strict 
global minimum for $G$ on $\ob$.
\qed

\medskip
\noindent
{\bf Proof of Theorem \EE.4.}
Let $H^t, \ t\geq0$ denote the family of downward-pointing $F$-harmonics at $x_0$ on $\O$
with $H^t\bigr|_{\bo} =\vf$, constructed in Theorem \EE.1.  Let 
$$
f(t) \ \equiv\ \T(H^t).
$$
To prove Theorem \EE.4 it is enough to show  that $f([0,\infty)) = [0,\infty)$,
and then choose $t\in f^{-1}(\T)$ and $H=H^t$.  The next lemma supplies this needed fact.

\def\hh{\ps}

\Lemma{\EE.6} 
{\sl
The function $f(t)$ satisfies:
\smallskip

(A) \ \ $f(0)=0$,
\smallskip

(B) \ \ $f(s) + (t-s) \leq f(t)$ for $0\leq s\leq t$,
\smallskip

(C) \ \ $f(t)$ is concave.
}
\pf
Part (A) is immediate from the assumption (H4) since $\Hz$ is $F$-harmonic on $\O$.
To prove (B) we show
\smallskip

(B)$'$ \ \ $H^s(x) - H^s(x_0) + (t-s) \ps(x) \ \leq\  H^t(x) - H^t(x_0)$.
\smallskip

\noindent
The functions
$$
u\ \equiv \ H^s_r(x) + (t-s) (\ps(x) -1) \and v\ \equiv\ H^t_r(x)
$$
have the same boundary values on $\partial B_r$ since 
$\hlb_s + (t-s) (\hh-1) = \hlb_t$, while on $\bo$ we have 
$u=\vf +(t-s)(\hh(x) -1) \leq \vf = v$.
By the convexity assumption on $F$, $u$ is $F$-subharmonic on $\O-B_r$.
Since $v$ is $F$-harmonic on  $\O-B_r$, comparison implies that $u\leq v$ on $\O-B_r$.
Taking $r\downarrow 0$ gives
$$
H^s(x)   + (t-s)( \hh(x)-1) \ \leq\  H^t(x)  
$$
which implies (B)$'$ since $H^s(x_0) = \Hz(x_0)-s$,  $H^t(x_0) = \Hz(x_0)-t$, and $\ps(x_0)=0$.

For (C) we show that for $0\leq \s \leq 1$ and $0<s<t$,
$$
(C)'\quad
\s \left( H^s(x) - H^s(x_0)  \right)  +  (1-\s)\left( H^t(x) - H^t(x_0)  \right)  
\ \ \leq\ \ 
 H^{\s s  + (1-\s)t}(x) - H^{\s s  + (1-\s)t}(x_0)  
$$
on $\ob$, which implies (C) by (H3).
Since $H^s(x_0) = H^0(x_0)-s$ for all $s$, the inequality (C)$'$ is equivalent to
$$
  \s H^s(x) + (1-\s) H^t(x) \ \leq \ H^{\s s + (1-\s)t} (x)\quad {\rm on} \ \ \ob.
\leqno{(C)''}
$$
This follows from 
$$
\s H^s_r(x) + (1-\s) H^t_r(x) \ \leq \ H^{\s s + (1-\s)t}_r (x)\quad {\rm on} \ \ \ob-B_r,
\eqno{(\EE.19)}
$$
which is true because the LHS and the RHS have the same boundary values on
$\partial(\O-\overline B_r)$ and, by the hypothesis that $F$ is a convex cone, 
the sum of the two $F$-harmonics on the LHS is $F$-subharmonic
and therefore in the Perron family for $H^{\s s + (1-\s)t}_r$.  \qed

\Remark{\EE.7}  Part  (4b) of Theorem \EE.1 implies $f(t) \leq t+\l(t)$
where $\l(t)$ is the affine function of $t$ used in the definition of $\hub_t$.
By  Lemma \EE.6 (B) we have $t\leq f(t)$.  We ask the question: When does $f(t)=t$?

\vskip.3in

\centerline
{\bf  Applications of Theorem \EE.4.
}

\medskip

There are many interesting examples of subequations $F$ and singularities $\ps$
for which the hypotheses (H3) and (H4) hold.

\medskip
\noindent
{\bf Case \EE.8. (Strong Uniqueness of Tangents).} 
Here we assume that  the  subequation $F$
 is invariant under ${\rm O}(n)$. Then $F$ has a well-defined Riesz
characteristic $p$ which we assume to satisfy $1\leq p<2$.
The function $\ps(x) = |x|^{2-p}$ provides the asymptotic singularity
Every $F$-subharmonic function $u$ defined near a point $x_0\in\rn$ has a well-defined
density $\T=\T(u,x_0)$ at $x_0$ defined to be $\lim_{r\to0} \sup_{B_r(x_0)} u / r^{2-p}$.
  We say that $F$ satisfies {\bf Strong Uniqueness of Tangents}
if, for each $u$,  every tangent to $u$ at $x_0$  is $\T(u, x_0)  |x-x_0|^{2-p}$.
This Strong Uniqueness holds for every $F$
 with the exception of $F=\cp$, which is the only possibility when $p=1$ (see   [\HLTangI], [\HLTangII]).

\medskip

Proposition 12.6  in [\HLTangI] states that for any $F$-subharmonic (not just a punctured 
$F$-harmonic)  function $u$, strong uniqueness of the tangent to $u$ holds if and only if 
(H3) holds, i.e., $u \sim \T K_p$.

This leads to a result where the several hypotheses are automatic

\Theorem{\EE.9} 
{\sl
Suppose $F$ is O$(n)$-invariant convex cone subequation with Riesz characteristic $1<p<2$
Take $\ps(x) = |x-x_0|^{2-p}$.
 Then   the (DPPD) is uniquely solvable for $F$.
}

\pf
The hypotheses (H3) and (H4) are true for $h$ so that Theorem \EE.4 applies.\qed

\medskip
\noindent
{\bf Case \EE.10. (Armstrong-Sirakov-Smart Fundamental Solutions)} 
Now we assume instead  that the  subequation $F$ in Theorem \EE.1 is uniformly elliptic,
and let $\ps=\Phi$ be the fundamental solution discussed in Example \BBB.5.

\Theorem {\EE.11. (Armstrong-Sirakov-Smart [AS$_1$])} {\sl
Every $F$-harmonic function $H$ on $B_r(x_0)-\{x_0\}$, which has a downward-pointing singularity at
$x_0$  is asymptotically equivalent to $\T \Phi(x-x_0)$ for some $\T\geq0$.
}
\medskip

This has the following corollary.

\Theorem{\EE.12} {\sl
Suppose the subequation $F$  is uniformly elliptic, and that the 
 downward-pointing fundamental solution $\Phi$ of Armstrong, Sirakov and Smart is of finite type
 (i.e.. of homogeneity $>0$). Then with $\ps=\Phi$ the hypotheses (H3) and (H4) are satisfied and Theorem \EE.4 applies.}


\vskip .3in
\centerline
{\bf \FF.  Prescribing Values at Singularities in the Finite Case.}
\medskip

\def\Th{ \gamma}

In the finite case one can also consider the  Dirichlet Problem with multiple prescribed singularities
where the value of the solution, instead of the density,  is given at each singular point. 
The $h_j$-density at each singular point $x_j$ is not precise but 
is replaced by  two-sided bounds on the difference quotient.
The general result is the following.
\medskip

\Theorem {\FF.1}  {\sl
Let $F$ be a cone subequation satisfying  (F3). Fix a function  $\vf\in C(\bo)$ 
and  points $x_1, ... , x_k \in \O$.   Suppose that  for each $j = 1, ... k$ we are given
a function
\smallskip
\centerline
{
$h_j\in C(\ob)$ which  is $F$-harmonic on $\O-\{x_j\}$
}
\smallskip
\noindent
with a downward-pointing singularity  of finite type at $x_j$.
Furthermore, assume that
\medskip

(i) \ \  \ the boundary of $\O$  is smooth and strictly  $F$-convex,

\medskip

(ii)\ \ \ each   $h_j$ has a strict global  minimum at $x_j$.

\medskip

(iii) \ there exists a function $h\in C(\ob)$ which is   $F$-subharmonic  on $\O-\{x_1,...,x_k\}$

\qquad\  with $h\sim h_j$ at $x_j$ for each $j$ and with $h\bigr|_{\bo} \leq \vf$.

\medskip
\noindent
 
\noindent
Then  there exists a function $H$ such that 
\medskip

(1) \quad $H\in C(\ob)$,
\medskip

(2)  \quad  $H$ is $F$-harmonic on $\O - \{x_1,...,x_k\}$,
\medskip

(3)  \quad  $H\bigr|_{\bo} = \vf$ \ \ and\ \ $H(x_j) = h(x_j)$\  for  $j=1, ... , k$,
\medskip

(4) \quad  There  exists a constant $c>1 $ such that for  any $\e>0$
$$
1-\e \ \leq\ {H(x) - h(x_j) \over h_j(x) - h_j(x_j)} \ \leq \ c\qquad {\rm for}\ x\ {\rm sufficiently\ near \ } x_j, \ \ j=1,...,k.
$$
}
\medskip
We postpone the proof   to the end of this section, and first examine some special cases.
For example, this theorem can be applied to subequations $F$ as in Case \EE.8 above
with the additional hypothesis of convexity.

\Cor{\FF.2}
{\sl
Suppose that $F$ is a convex, O$(n)$-invariant subequation whose Riesz characteristic 
$p$ satisfies $1<p<2$.  Let $\O\ss\ss\rn$ be a domain with a smooth strictly $F$-convex boundary.
Then given points $x_1, ... , x_k \in \O$, positive numbers $\Th_1, ... , \Th_k >0$,
and a function $\vf\in C(\O)$, the following holds.  For every constant $C$ such that the
restriction
$$
\left.
\left\{h\ \equiv\ h_C \ \equiv\ \sum_{j=1}^k \Th_j |x-x_j|^{2-p} +C \right\} \right|_{\bo} \ \leq\ \vf,
\eqno{(\FF.1)}
$$
there exists a function $H$ with properties (1) -- (4) above.
}
\medskip
\noindent
{\bf Proof.}   Since $F$ is convex, the function $h$ is  $F$-subharmonic.  It satisfies the 
boundary hypothesis in (iii) by (\FF.1). 
Finally, one has that $h(x)\sim \Th_j |x-x_j|^{2-p}$ at $x_j$  for each $j$. 
(To see this note that if $h(x) = \T |x|^\a + g(x)$
 where $g$ is smooth and $0<\a<1$, then $(h(x)-h(0))/|x|^\a \to \Theta$ as $x\to 0$.) 
 Hence, Theorem \FF.1 applies.\qed
\medskip

 The simplest subequation of this type is $\cp_p$.

We note that in the case where $p=1$  ( i.e., $F=\cp$), multiple singularities
cannot exist, because a convex function cannot have more than one
strict local minimum.  Thus Theorem \FF.1 does not apply since the global function
$h$ does not exist when $k>1$. 

 Corollary \FF.2 allows us to prescribe certain values for the function $H$ at the singular points.
 This can be thought of as the Dirichlet problem, 
if one considers $x_1, ... , x_k$ as additional boundary points.
Uniqueness holds for this problem by comparison (see [\HLAE, Thm. 6.2]) , but 
in general  the values at these interior  points cannot be given  arbitrarily (see Proposition \FF.5 below).
However, Corollary \FF.2 does provides a set $\cv\ss\bbr^k$ of values which can 
be prescribed. More precisely,
let $F$, $\O$, $\vf$ and $x_1, ... , x_k$ be as in Corollary \FF.2.
Given  $\Th_1, ... , \Th_k> 0$ and a constant $C$, let 
$$
h_{\Th,   C}(x) 
\equiv \sum_j \Th_j |x-x_j|^{2-p} +C.
$$ 
Define  the  set  $\cv\ss \bbr^k$ by the condition that $=(v_1, ... , v_k) \in\cv$ iff
$$
v_j \ =\ h_{\Th, C}(x_j)  \quad{\rm for\ some\ }\ \Th, C \ \ {\rm such\ that\ \ } h_{\Th, C}\bigr|_{\bo}\leq \vf.
\eqno{(\FF.2)}
$$

\Cor{\FF.3.  ($\cv\ss  \cv \! {\it al}$)}
{\sl  For every $v\in \cv$ there exists a solution $H$ to the ``value problem''.  That is, there exists
 $H\in C(\ob)$  which is    $F$-harmonic on $\O-\{x_1,...,x_k\}$, and takes on the boundary values
 $$
H\bigr|_{\bo} \ =\ \vf\and H(x_j) \ =\ v_j \ \ \ {\rm for}\ \ j=1,...,k.
 $$
In addition, $H$ satisfies (4) above  with $h_j(x) = \Th_j|x-x_j|^{2-p}$.
Moreover the set $\cv$ has non-empty interior and satisfies $ \cv +(c,...,c) \ \ss\ \cv$ for all $c\leq0$.}
 
 \medskip
\noindent
{\bf Proof.} If $v$ is given as in (\FF.2) with
$h_{\Th, C}\bigr|_{\bo} < \vf$, then $h_{\Th', C'}\bigr|_{\bo} < \vf$
for $(\Th',C')$ in a neighborhood of $(\Th,C)$ and the resulting values $v'$ 
fill out a neighborhood of $v$ in $\bbr^k$.  To see this consider the symmetric
$k\times k$ matrix $A = ((a_{ij}))$ where $a_{ij} = |x_i-x_j|^{p-2}$. Then $\det(A)\neq0$ and
so the mapping $\Th \mapsto A\cdot \Th$ is open.  This proves that 
$\cv$ has non-empty interior. The remaining assertions are clear.\qed

\Note{\FF.4. (The Nonlinear Green's Functions)} 
In the case of a single  point $x_0\in\O$ and boundary values $\vf =  0$,
Corollary  \FF.3 produces
the family of Green's functions $G(x) = G_\O(x; x_0,\Theta |x-x_0|^{2-p})$
in  Corollary \EE.5. To see this, note that the function $H(x)$ given
in Corollary \FF.3 is $\leq 0$ on $\ob$ since by (4) it has no test functions
at $x_0$ and is therefore $F$-subharmonic on $\O$. If $H(x_0)=0$, then $H\equiv 0$ 
on $\ob$ by the maximum principle on $\O-\{x_0\}$ applied to the $\ft$-subharmonic
function $-H$. Thus $H(x_0)<0$ and by uniqueness $H(x) = G(x;x_0, \Theta |x-x_0|^{2-p})$
where $\Theta\geq \Th >0$ is the density of $H$ at $x_0$.
It  is now obvious by rescaling (since  $\vf=0$) that the set $\cv$ defined by (\FF.2) at $x_0$ is exactly
$\{v\leq0\}$.
\medskip

There remains the question of describing $\cv \! {\it al}$, i.e., when, if at all, can one prescribe $v=H(x_0)>0$.

\Prop{\FF.5} 
{\sl
Suppose $H\in C(\ob)$ is $F$ harmonic on $\O-\{x_0\}$ and satisfies $H\bigr|_{\bo}=0$.
Then $H(x_0)\leq0$.  Thus,  $\cv=\{v\leq 0\}$ is exactly the set of values $\cv \! {\it al}$ that can be prescribed
for the one-point Dirichlet problem with $\vf=0$.
}
\pf
Assume $x_0=0$ for simplicity, and suppose that $H(0)>0$.
Choose $R>0$ so that $\ob\ss\ss B_R(0)$, and  for $q>2$ and $r>0$
consider the  function
$$
v_r(x)  \  \equiv \ r \left \{ {1\over |x|^{q-2}} -{1\over R^{q-2}}  \right\}.
$$
The ordered eigenvalues of $D^2_x (v_r)$, up to the positive factor ${r(q-2)\over |x|^q}$, are  
$$
- 1, \ ...\ , - 1, (q-1)
 $$
and so the eigenvalues of  $D^2_x (v_r)$ satisfy
$$
\l_{\rm min} +(p-1)\l_{\rm max} \ \cong \ -1 + (p-1)(q-1) \ <\ 0
\eqno{(\FF.3)}
$$
if we choose
$$
q\ <\  1+{1\over p-1}
$$
which is possible since $0< p-1<1$.

Now for $r>0$ sufficiently large we have $v_r > H$ on $\ob$.
(In fact by the maximum principle on $\O-\{0\}$ one has that $H(0)$ is a global maximum
of $H$ on $\ob$.)  Let
$$
r_0 \ \equiv \ \inf \{r : v_r > H \ \ {\rm on\ \ } \ob\}\ \geq\ 0.
$$
Now $r_0$ cannot be 0 since we are assuming that $H$ is continuous with
$H(0)>0$.  Thus, there is a point $y\in \ob$ where $H(y)=v_{r_0}(y)$.
Note that  $y\notin \bo$  since $H=0$ there,  and $y\neq0$ since $v_{r_0}(0)=\infty$.

We conclude that  $y\in\O-\{0\}$ and $v_{r_0}$ is a test function for $H$ at $y$. In particular,
$D^2_y(v_{r_0}) \in F$. 
However,  since $F$ is  O$(n)$-invariant and of characteristic $p$,  we have
by  Proposition 3.13 in [\HLTangI] that  $F\ss \cp_p^{\rm min/max}
\equiv \{ \l_{\rm min} +(p-1)\l_{\rm max} \geq 0\}$.
Thus we have $D^2_y (v_{r_0})  \in  \cp_p^{\rm min/max}$ contradicting (\FF.3).\qed

\medskip

For multiple singular points things are much more complicated.

\Cor{\FF.6. (Prescribing Values at Multiple Singularities)} 
{\sl
Let $F$ and $\O$ be as above, and  fix  points  $x_1,...,x_k\in\O$.
Take  boundary values $\vf=0$ and choose $\Th \in \bbr^k_+$. 
Then for each $C$ satisfying  $h_{\Th, C}\bigr|_{\bo}\leq0$,
 one obtains from
Corollary \FF.3 the existence of a   {\bf multi-pole  Green's function} $G(x) = G_\O(x; x_1,...,x_k)$ 
with a singularity  of type $|x-x_j|^{2-p}$ at $ x_j$ but with an unknown density $\Theta_j\geq \Th_j$ at $x_j$.
As in (\FF.2)  let  $\cv\ss\bbr^k$  be the set of all  values 
$v= (G(x_1),..., G(x_k))$ at the singular  points, obtained by varying $\Th$.
 Then $\cv$ is a convex cone  in the negative ``octant'' $\bbr_-^k\ss\bbr^k$.}

\pf
The set $\cv$  given in Corollary \FF.3 is described as follows.  For each 
$\Th = (\Th_1,...,\Th_k) \in \bbr^k_+$ set 
$$
 F_{\Th }(x) \ =\ \sum_{j=1}^k \Th_j |x-x_j|^{2-p}
 \and C(\Th) \ =\ - \sup_{x\in \bo} F_\Th(x)
$$
Then $v=(v_1,...,v_k)\in \cv$ if and only if there exist $\Th$ and $C$ such that
$$
v_j \ =\  F_\Th(x_j) - C \quad{\rm for \ \ } j=1,..., k \and C\leq C(\O).
\eqno{(\FF.5)}
$$
 It is clear that $C(t\O) = t C(\O)$ and $F_{t\Th}= tF_\Th$ for all $t\geq0$, so $\cv$
 is a cone with vertex the origin.  
It remains  to show that $\cv+\cv\ss\cv$. We begin by observing that 
 $$
 C(\Th + \Th') \ \geq \ C(\Th) + C(\Th')
 \eqno{(\FF.6)}
$$
 (since the sup of the sum is $\leq$ the sum of the sups).
 Consider the vectors $v(\Th)$ obtained by taking $C=C(\Omega)$ in (\FF.5).
 Then  we have
 $$
 \eqalign
 {
 v_j(\Th + \Th') \ &=\ F_{\Th + \Th'}(x_j) + C(\Th + \Th')   \cr
 &= \ F_{\Th}(x_j) + F_{\Th'}(x_j) + C(\Th) + C(\Th') -\kappa  \cr
 \ &=\ v_j(\Th)+v_j(\Th') -\kappa
 }
 $$
 where $\kappa \equiv  C(\Th + \Th') -C(\Th) + C(\Th') \geq0$ by (\FF.6).
 Thus by (\FF.5) we have $v(\Th) + v(\Th') \in \cv$, and the assertion follows easily.
 \qed
 \medskip
 
For a slightly different perspective, note that $F_\Th(x) = \Th\cdot F(x)$ where
 $$
 F(x) = (|x-x_1|^{2-p}, |x-x_2|^{2-p}, ...  , |x-x_k|^{2-p}).
 $$
 and consider the symmetric $k\times k$ matrix with positive entries:
 $$
A\ =\ \left( \matrix{F(x_1) \cr F(x_2) \cr : \cr F(x_k)}\right) \ =\ (F(x_1)^t, F(x_2)^t, ... ,  F(x_k)^t)
$$
Then the set $\cv$ is given by:
$$
v\in\cv
\qquad\iff\qquad 
v\ =\ A\cdot \Th + C(1,...,1) \quad{\rm for }\ \ \Th \in \bbr_+^k \ \ {\rm and}\ \ C\leq C(\O)
\eqno{(\FF.7)}
$$

Taking $\Th = e_j = (0,...,0,1,0,...,0)$ (the $j^{\rm th}$ coordinate vector in $\bbr^k$)
for $j=1,...,k$ we have the following.

\Lemma {\FF.7} 
 {\sl
The set $\cv$ contains the convex cone in the negative octant generated by the vectors
$$
V_j \ \equiv\ F(x_j) -\sup_{x \in \partial \O} |x-x_j|^{2-p} (1,1, ..., 1), \qquad  j=1,...,k.
$$
}

\Remark {\FF.8}
It is important to note that the functions $G_\O(x;x_1,...,x_k)$ constructed in 
Corollary \FF.6 do not have $h_j$- density equal to $\Th_j$ at $x_j$ (where
$h_j(x) = |x-x_j|^{2-p}$).  By Theorem \FF.1 (4) the constant $\Th_j$ is a lower
bound on the density and there exists a global upper bound depending on the
data.

\medskip
\noindent
{\bf Question \FF.9.} Can one determine these densities, at least in geometrically
simple cases?

\medskip
\noindent
{\bf Question \FF.10.} Can one determine  the value  set $\cv$ defined by (\FF.2), at least in relatively
simple cases?

\medskip
\noindent
{\bf Question \FF.11.} What is the full set $\cv\!{\it al}  \supset \cv$ of possible values $v$ (both depending
on $\vf \in C(\bo)$) for which there exists a solution $H$ to the ``value problem'' (as defined in Corollary
\FF.3).

\medskip
 Note \FF.4 and Proposition \FF.5 answer this question  in the case of a single point.

\medskip
\noindent
{\bf Proof of Theorem \FF.1.} We apply the existence construction given in Section 6.
Consider $r>0$ sufficiently small that the closed  balls $\overline B_r(x_j), 1\leq  j\leq k$
are mutually disjoint and contained in $\O$. Let $H_r$ be the $F$-harmonic function
 on $\O_r \equiv \O - \bigcup_j \overline B_r(x_j)$ with boundary values $\vf$ on $\bo$
 and $h$ on each $\partial  B_r(x_j)$. We extend $H_r$ to $\overline H_r\in C(\ob)$ 
 by setting 
 $\overline   H_r = h$  on each $B_r(x_j)$. 
 The arguments of Section  \AA \ show that 
 $$
\overline  H_r \ \ \uparrow\ \ H \qquad{\rm as}\ \ r\to 0.
 $$
As before $H^*$ is $F$-subharmonic and $-H$ is $\ft$-subharmonic on 
$\O-\{x_1,...,x_k\}$. As before, let $H^{\rm DP}$ be the solution to the
Dirichlet problem on $\O$ with  boundary values $\vf$.  Then $\overline   H_r\leq  H^{\rm DP}$ 
and therefore $H\leq H^{\rm DP}$ on $\ob$. Thus $H^*\leq H^{\rm DP}$ and in particular $H^*\bigr|_{\bo} =\vf$.
Note also that $h\leq H\leq H^*$.

Now for each $j$ we choose a constant $c_j >0$  sufficiently large that the punctured $F$-harmonic function
$$
H^j(x) \ \equiv\ c_j(h_j(x)  - h_j(x_j)) + h(x_j)
$$
satisfies
$$
H^j \ >\ h \quad{\rm on}\ \ \ob - \{x_j\}
\and
H^j \ >\ \vf \quad{\rm on}\ \ \bo.
$$
This can be done since $h_j$ has a strict global minimum at $x_j$.
It follows that $\overline H_r \leq H^j$ on $\ob$ for all small $r$, and therefore $H\leq H^j$. 
In sum we have that
$$
h\ \leq\ H\ \leq\ H^*\ \leq\ H^j\quad{\rm for\ each}\ j=1,...,k,
$$
and in particular $H$ and $H^*$ are continuous at each $x_j$ with 
$H(x_j)=H^*(x_j) = h(x_j)$.  Subtracting $h(x_j)$ and dividing by the positive function $h(x)-h(x_j)$ gives 
$$
1 \ \leq {H(x)-h(x_j) \over h(x)-h(x_j)} 
\ \leq\ {H^*(x)-h(x_j) \over h(x)-h(x_j)} 
\ \leq\ {H^j(x)-h(x_j) \over h(x)-h(x_j)} \ =\ c_j \ \ {\rm near}\  x_j
\eqno{(\FF.8)}
$$

We now consider the function $u \equiv H^*-H$ which is continuous on $\ob$, subaffine on $\O-\{x_1, ... , x_k\}$,
and satisfies $u\bigr|_{\bo} = 0$ and $u(x_j)=0$ for all $j$. It follows that $u\leq 0$ on $\ob$. 
Hence, $H=H^*$ and so $H$ is $F$-harmonic on $\O-\{x_1, ... , x_k\}$.

The assertion (4) now follows from (\FF.8) and the fact that $h\sim h_j$ at $x_j$, i.e., 
$$
\lim_{x\to x_j} {h(x)-h(x_j) \over h_j(x)- h_j(x_j)} \ =\ 1. \qquad\mathqed
$$

\vskip.3in


\centerline{\headfont  Appendix A.   Asymptotic Equivalences and Tangent Flows (for p $\ne$ 2).}
\bigskip

The asymptotic equivalence classes can be related to the tangent flow very generally.
Suppose $u$ is an upper semi-continuous  function defined in a neighborhood
of the origin. Recall the tangent flow and the Riesz kernel:
$$
u_r(x)\  \equiv \ 
\cases
{
r^{p-2} u(rx) \qquad {\rm if}\ \ p\geq  2  \cr
{u(rx)-u(0)  \over r^{2-p}} \qquad {\rm if}\ \ 1\leq p< 2.
}
\and
K(x) \ \equiv\ 
\cases
{
-{1\over |x|^{p-2}} \qquad {\rm if}\ \ p\geq  2  \cr
|x|^{2-p} \qquad {\rm if}\ \ 1\leq p< 2.
}
$$
Note that $K_r = K$ for the full range $p\geq 1, \ p\neq 2$.

Fix $\Theta >0$.
We define  the first notion of asymptotic equivalence at 0 as follows.

\Def{A.1}  
$\qquad u\ \sim\ \T K$ at 0 \ \ if $$
\eqalign
{
&\ \ \ \lim_{x\to0} {u(x)\over K(x) }\ =\ \T \qquad \qquad {\rm when}\ \ p\geq 2  \cr
&\lim_{x\to0} {u(x)-u(0)\over K(x) }\ =\ \T\qquad \   {\rm when}\ \ 2>p\geq 1  \cr
}
\eqno{(A.1)}
$$

 We have a second alternative definition of asymptotic equivalence when $p>2$.
 
 \Def{A.2} $\qquad u\ \approx\ \T K$ at 0 \ \ if 
$$
u-\T K \quad {\rm is\ bounded\ in\ a\ neighborhood\ of \ the\ origin.}
\eqno{(A.2)}
$$
This notion is stronger.  (When $2 >p\geq 1$ it is weaker, in fact too weak to be useful.)

\Prop{A.3. ($p>2$)}
{\sl If $u\approx \T K$, then $u\sim \T K$. However, the converse is false.}
\pf
Note that
$$
\eqalign
{
u\ \approx\ \T K 
\qquad &\iff\qquad
|u(x)-\T K(x)| \ \leq\  C \ \ {\rm near}\ x=0   \cr
\qquad &\iff\qquad
\left| {u(x) \over K(x)} -\T \right| \ \leq {C \over |K(x)|} \ \ {\rm near}\ x=0
}
$$
This implies that  $\lim_{x\to0} u(x)/K(x)=\T$ since $K(0)=-\infty$.  See Example A.6 for
a counterexample to the converse.\qed
\medskip

These asymptotic equivalences are related to the tangent flows as follows.

\Prop{A.4. ($p>2$)}
$$
\eqalign
{
& (a)\ \ \ u\ \approx \ \T K\qquad \Rightarrow\qquad u_r\  \to \  \T K\ \ {\rm uniformly\ on\ \ } \overline B_R \cr
& (b)\ \ \ u_r -  \T K \ \ {\rm bounded\ near\  }\ \ 0\ \ {\rm for \ some\ } r
\qquad\Rightarrow \qquad u\ \approx\ \T K.
}
$$

\Prop{A.5}
$$
\eqalign
{
& (a)\ \ (1\leq p<2) \qquad  u\ \sim \ \T K\qquad 
\Rightarrow\qquad u_r\  \to \  \T K\ \ {\rm uniformly\ on\ \ } \overline B_R \cr
& (b)\ \ (2<p)\qquad \qquad u\ \sim \ \T K\qquad 
\Rightarrow\qquad u_r\  \to \  \T K\ \ {\rm uniformly\ on\ \ } \overline A_{s,R} \cr
& (c)\ \ (1\leq p<\infty)\qquad 
u_r\  \to \  \T K\ \ {\rm uniformly\ on\  some\ sphere\ } \partial B_R
\qquad  \Rightarrow\qquad 
 u\ \sim \ \T K \cr
}
$$

\medskip
\noindent
{\bf Proof of Proposition A.4(a).}  
Assume $u\approx \T K$.  The inequality
$$
|u(y) - \T K(y) |\ \leq \ C
\eqno{(A.3)}
$$
can be rewritten with $y=rx$ as
$$
|u_r(x) - \T K(x) |\ \leq \ Cr^{p-2}
\eqno{(A.3)'}
$$
by applying the tangent flow to both sides.  Thus if
(A.3) holds for all $|y|\leq\d$, then (A.3)$'$ holds for all $|x|\leq R$ and $r\leq \d/R$,
which suffices to prove that $u_r$ converges uniformly to $\T K$ on $B_R$.

\medskip
\noindent
{\bf Proof of Proposition A.5(a), (b).}  
With $y,r,x$ related by $y=rx$, the inequalities
$$
(a)\ \ \ 
\left | {u(y) \over K(y)} - \T \right|\ \leq \ \e
\and
(b)\ \ \ 
|u_r(x) - \T K(x) |\ \leq \ \e |K(x)|
\eqno{(A.4)}
$$
are equivalent.
If (a) holds for $|y|\leq \d$, then (b) holds for $|x|\leq R$ and $r\leq \d/R$.

\medskip
\noindent
{\bf Case: $1\leq p<2$.}
Note that $|x|\leq R \Rightarrow |K(x)| = K(x) \leq K(R)$.
This proves that (A.4a) for $|y| \leq \d$ implies $|u_r(x)-\T K(x)|\leq \e K(R)$ 
for $|x|\leq R$ and $r<\d/R$.  Thus $u_r$ converges uniformly to $\T K$ on $\overline B_R$.

\medskip
\noindent
{\bf Case: $2<p$.}
Note that $s \leq |x| \Rightarrow |K(x)| \leq |K(s)|$.
This proves that (A.4a) for $|y| \leq \d$ implies that 
 $|u_r(x)-\T K(x)|\leq \e |K(s)|$ 
for $s\leq |x|\leq R$ and $r<\d/R$, 
since $|y|=r|x|\leq {\d\over R} R =\d$.
 Thus $u_r$ converges uniformly to $\T K$ on $A_{s,R}$.

\medskip
\noindent
{\bf Proof of Proposition A.4(b).}  
With $y,r,x$ related by $y=rx$, the inequalities
$$
\eqalign
{
&(a)\ \ \  \ 
|u_\d(x) -\T K(x)|\ \leq\ C\qquad\quad \forall \ |x|\leq R
\and\cr
&(b)\ \ \ 
\d^{p-2}|u(y) -\T K(y)|\ \leq\ C\qquad\forall \ |y|\leq  \d R
}
\eqno{(A.5)}
$$
are equivalent.
If (a) holds for all $|x|\leq R$, then (b) holds for all $|y|\leq \d R$.
This  proves that $u \approx \T K$ if (a) is true for some $\d, C>0$ and all $|x|\leq R$

\medskip
\noindent
{\bf Proof of Proposition A.5(c).}  
Again if $y=rx$, then
$$
(a)\ \ \ 
|u_r(x) - \T K(x) |\ \leq \ \e
\and
(b)\ \ \ 
\left | {u(y) \over K(y)} - \T  \right |\ \leq \ {\e \over |K(x)|}
\eqno{(A.6)}
$$
are equivalent.  To see this divide both sides of (A.6a) by $|K(x)|$ and note that 
$
{u_r(x) \over K(x)} ={u_r(x) \over K_r(x)} = {u(rx) \over K(rx)}.
$
Suppose (A.6a) holds for all $|x| = R$ and $r\leq \d$.
Then (A.6b) holds for all $y=rx$ with $|x|=R$ and $r\leq\d$,
or for all $|y|\leq \d R$.  Since $K(x)= K(R)$ this is enough to prove
that $\lim_{y\to0} u(y)/K(y) =\T$.

\Ex{A.6. (Counterexamples to:\  $u \sim \T K \  \Rightarrow\  u \approx \T K$)}
Take 
$$
F=\D, \qquad  K(x) = -{1\over |x|^{n-2}},  \and u(x) \equiv -{\T\over |x|^{n-2}} - {1\over |x|^{p-2}} 
$$
with $p< n$. Then
$
\lim_{x\to0} {u(x)\over K(x)} \ =\ \lim_{x\to0}\left( \T + |x|^{n-p}   \right) \ =\ \T
$
so that $u \sim K$.  However, $u-\T K = -{1\over |x|^{p-2}}$ 
is not bounded near the origin, i.e., 
$u\approx \T K$ is not true.  Note that $u$ is $\D$-subharmonic.

In this example  we could also replace $n$ with any $q\leq n$ and take $p<q$.  In this case
the corresponding function $u$ is $\cp_q$-subharmonic since $\cp_p\ss\cp_q$.

\vskip.4in


\centerline{\bf REFERENCES}

\vskip .2in

\noindent
\item{[A$_1$]}   S. Alesker,  {\sl  Non-commutative linear algebra and  plurisubharmonic functions  of quaternionic variables}, Bull.  Sci.  Math., {\bf 127} (2003), 1-35. also ArXiv:math.CV/0104209.  

\smallskip

\noindent
\item{[A$_2$]}   \ \----------,   {\sl  Quaternionic Monge-Amp\`ere equations}, 
J. Geom. Anal., {\bf 13} (2003),  205-238.

 ArXiv:math.CV/0208805.  

\smallskip

\noindent
\item{[AV]}    S. Alesker and M. Verbitsky,  {\sl  Plurisubharmonic functions  on hypercomplex manifolds and HKT-geometry},  J. Geom. Anal. {\bf 16} (2006), no. 3, 375Ð399.

\smallskip

\noindent
\item{[AS$_1$]}  S. N.   Armstrong, B.  Sirakov  and  C. K. Smart,  {\sl Fundamental
solutions of homogeneous fully nonlinear elliptic equations}, Comm.
Pure. Appl. Math., 64 (2011), no. 6, 737-777.

\smallskip

\noindent
\item{[AS$_2$]}    \ \----------,   {\sl Singular solutions
of fully nonlinear elliptic equations and applications}, Arch. Ration.
Mech. Anal., 205 (2012), no. 2, 345-394.

\smallskip

\noindent
\item{[BT]}   E. Bedford and B. A. Taylor,  {\sl The Dirichlet problem for a complex Monge-Amp\`ere equation}, 
Inventiones Math. {\bf 37} (1976), no.1, 1-44.

\smallskip

\noindent
 \item{[CC]}    L. Caffarelli and X. Cabr\'e,  {Fully Nonlinear Elliptic Equations}, 
{\sl Colloquium Publications}, {\bf 43}, American Math. Soc., 1995.

\smallskip

\noindent
 \item{[CNS]}    L. Caffarelli, L. Nirenberg and J. Spruck,  {\sl
The Dirichlet problem for nonlinear second order elliptic equations, III: 
Functions of the eigenvalues of the Hessian},  Acta Math.
  {\bf 155} (1985),   261-301.

 \smallskip

\noindent
\item{[C]}   M. G. Crandall,  {\sl  Viscosity solutions: a primer},  
pp. 1-43 in ``Viscosity Solutions and Applications''  Ed.'s Dolcetta and Lions, 
SLNM {\bf 1660}, Springer Press, New York, 1997.

 \smallskip

\noindent
\item{[CIL]}   M. G. Crandall, H. Ishii and P. L. Lions, {\sl
User's guide to viscosity solutions of second order partial differential equations},  
Bull. Amer. Math. Soc. (N. S.) {\bf 27} (1992), 1-67.

 \smallskip

\noindent
\item{[GLN]}   B. Gidas, W. M. Li and L. Nirenberg,   {\sl
Symmetry and related properties via the maximum principle},  
Comm.  Math.  Phys.   {\bf 68} (1979), 209-243.

 \smallskip

\item {[\HLDD]}   \ \----------,  {\sl  Dirichlet duality and the non-linear Dirichlet problem},    Comm. on Pure and Applied Math. {\bf 62} (2009), 396-443. ArXiv:math.0710.3991

\smallskip

\item {[\HLPUP]}  \ \----------,    {\sl  Plurisubharmonicity in a general geometric context},  Geometry and Analysis {\bf 1} (2010), 363-401. ArXiv:0804.1316.

\smallskip

\item {[\HLDDR]}  \ \----------,   {\sl  Dirichlet duality and the nonlinear Dirichlet problem
on Riemannian manifolds},  J. Diff. Geom. {\bf 88} (2011), 395-482.   ArXiv:0912.5220.
\smallskip


\item {[\HLREST]}  \ \----------, {\sl  The restriction theorem for fully nonlinear subequations}, 
   Ann. Inst.  Fourier {\bf 64}  No. 1 (2014) , 217-265. 
ArXiv:1101.4850.

\smallskip

\item {[\HLPCON]}   \ \----------,  {\sl  p-convexity, p-plurisubharmonicity  and the Levi problem },
   Indiana Univ. Math. J.  {\bf 62} No.\ 1 (2013), 149-169.  ArXiv:1111.3895.

\smallskip

\item {[\HLSURVEY]}   \ \----------,  {\sl  Existence, uniqueness and removable singularities
for nonlinear partial differential equations in geometry},  pp. 102-156 in ``Surveys in Differential Geometry 2013'', vol. 18,  
H.-D. Cao and S.-T. Yau eds., International Press, Somerville, MA, 2013.
ArXiv:1303.1117.

 \smallskip

\item  {[\HLRS]} \ \----------, {\sl  Removable singularities for nonlinear subequations}, Indiana Univ. Math. J., {\bf 63},
No. 5 (2014), 1525-1552.
 ArXiv:1303.0437.
\smallskip

\item  {[\HLBP]} \ \----------, {\sl The equivalence of  viscosity and distributional
subsolutions for convex subequations -- the strong Bellman principle},
Bull. Braz. Math. Soc.  (N.S.) {\bf 44} No. 4 (2013),  621-652.   ArXiv:1301.4914.

\smallskip

\noindent
\item  {[\HLAE]} \ \----------, {\sl The AE Theorem and Addition Theorems for quasi-convex functions}, 
 ArXiv: 1309:1770.

\smallskip

\noindent
\item  {[\HLTangI]} \ \----------,{\sl  Tangents to subsolutions -- existence and uniqueness, Part I},  ArXiv:1408.5797.

\smallskip

\noindent
\item  {[\HLTangII]} \ \----------, {\sl  Tangents to subsolutions -- existence and uniqueness, Part II},  ArXiv:1408.5851.

\smallskip

   \noindent
\item{[I]}    H. Ishii,  {\sl On the equivalence of 
two notions of weak solutions, 
viscosity solutions and distribution solutions},
Funkcial. Ekvac. 38 (1995), no. 1, 101Ð120.

   \smallskip

 \noindent
\item{[K$_1$]} 
M. Klimek, 
Pluripotential theory.
London Mathematical Society Monographs. New Series, 6. Oxford Science Publications. The Clarendon Press, Oxford University Press, New York, 1991.

   \smallskip

 \noindent
\item{[K$_2$]} 
 \ \----------, {\sl
Extremal plurisubharmonic functions and invariant pseudodistances}, 
Bull. Soc.  math. de France, {\bf 113} (1985), 231-240.

\smallskip

\item{[La$_1$]}  D.Labutin,
{\sl Isolated singularities for fully nonlinear elliptic equations},  J. Differential Equations  {\bf 177} (2001), No. 1, 49-76.

 \smallskip
 
\item{[La$_2$]}   \ \----------, {\sl Singularities of viscosity solutions of fully nonlinear elliptic equations}, 
Viscosity Solutions of Differential Equations and Related Topics, Ishii ed., RIMS K\^oky\^uroku
No. 1287, Kyoto University, Kyoto (2002), 45-57

\smallskip

\item{[La$_3$]}   \ \----------, {\sl Potential estimates for a class of fully nonlinear elliptic equations}, 
Duke Math. J. {\bf 111} No. 1 (2002), 1-49.

\smallskip

\item {[L]}   N. S.  Landkof,   {Foundations of Modern Potential Theory},  Springer-Verlag, New York, 1972.

\smallskip

\item{[Le]}  P. Lelong,
{\sl Fonction de Green pluricomplexe et lemmes de Schwarz dans les espaces de Banach},   J. Math. Pures Appl. (9) 68 (1989), no. 3, 319Ð347.

\smallskip

\item{[Lem]}  L. Lempert,
{\sl Solving the degenerate Monge-Amp\`ere equation with
one concentrated singularity},  Mathematische Annalen {\bf 263}  (1983), 515-532.

\smallskip

\item {[Sh]} J.-P. Sha, {\sl  $p$-convex riemannian manifolds},
Invent.  Math.  {\bf 83} (1986), 437-447.

\smallskip

   \noindent
\item{[TW$_1$]} 
N.  Trudinger and X-J.  Wang,  {\sl Hessian measures. I},  Dedicated to Olga Ladyzhenskaya. Topol. Methods Nonlinear Anal. {\bf 10} (1997), no. 2, 225--239.

\smallskip

   \noindent
\item{[TW$_2$]} 
  \ \----------,   {\sl Hessian measures. II},  Ann. of Math. (2)  {\bf 150} (1999), no. 2, 579--604.

\smallskip

   \noindent
\item{[TW$_3$]} 
  \ \----------,   {\sl Hessian measures. III},   J. Funct. Anal. {\bf  193} (2002), no. 1, 1--23.

\smallskip

\item {[Wu]}   H. Wu,  {\sl  Manifolds of partially positive curvature},
Indiana Univ. Math. J. {\bf 36} No. 3 (1987), 525-548.
 
  \smallskip

\item {[Z]}   A. Zeriahi,  {\sl  Pluricomplex Green functions and the Dirichlet problem for the complex
Monge-Amp\`ere operator},
Michigan Math. J. {\bf 44}  (1997), 579-596.

\vfill\eject

\end